\documentclass[11pt,a4paper]{amsart}

\usepackage{amssymb}
\usepackage{amscd}
\usepackage{fancyhdr}
\input xy 
\xyoption{all}

\usepackage{amsmath,amsfonts,stmaryrd}
\usepackage[dvipsnames,svgnames]{xcolor}
\usepackage[colorlinks, urlcolor=Navy, linkcolor=Navy, citecolor=Navy]{hyperref}
\usepackage[abbrev,msc-links]{amsrefs}
\usepackage[latin1]{inputenc}
\usepackage{tikz}
\usepackage{tikz-cd}
\usetikzlibrary{shapes,arrows}

\theoremstyle{plain}
\newtheorem{thm}{Theorem}[section]
\newtheorem{lem}[thm]{Lemma}

\newtheorem{pro}[thm]{Proposition}
\newtheorem{cor}[thm]{Corollary}

\theoremstyle{definition}
\newtheorem{dfn}[thm]{Definition}
\newtheorem{exa}[thm]{Example}
\newtheorem{rem}[thm]{Remark}

\DeclareMathOperator{\Hom}{Hom}

\DeclareMathOperator{\End}{End}
\DeclareMathOperator{\Ext}{Ext}

\DeclareMathOperator{\modu}{mod}

\DeclareMathOperator{\Modu}{Mod}

\DeclareMathOperator{\eff}{eff}

\DeclareMathOperator{\id}{id}
\DeclareMathOperator{\idim}{id}

\DeclareMathOperator{\gldim}{gldim}

\newcommand{\oTo}{\xymatrix{ \ar@{^{(}->}[r]|{\mathbf{O}}& }} 
\newcommand{\cTo}{\xymatrix{ \ar@{^{(}->}[r]|{\mathbf{|}}& }} 
\newcommand{\coTo}{\xymatrix{ \ar@{^{(}->}[r]|{\mathbf{O}}|{\mathbf{|}}& }} 

\DeclareMathOperator{\Ker}{ker}
\DeclareMathOperator{\coKer}{coker}
\DeclareMathOperator*{\colim}{colim}
\DeclareMathOperator{\Bild}{Im}

\newcommand{\defl}{\twoheadrightarrow}
\newcommand{\infl}{\rightarrowtail}

\newcommand{\La}{\Lambda}

\newcommand{\mcA}{\mathcal{A}}
\newcommand{\mcB}{\mathcal{B}}
\newcommand{\mcC}{\mathcal{C}}
\newcommand{\mcD}{\mathcal{D}}
\newcommand{\mcE}{\mathcal{E}}
\newcommand{\mcF}{\mathcal{F}}

\newcommand{\mcH}{\mathcal{H}}
\newcommand{\mcI}{\mathcal{I}}

\newcommand{\mcM}{\mathcal{M}}

\newcommand{\mcP}{\mathcal{P}}
\newcommand{\mcQ}{\mathcal{Q}}
\newcommand{\mcR}{\mathcal{R}}
\newcommand{\mcS}{\mathcal{S}}
\newcommand{\mcT}{\mathcal{T}}
\newcommand{\mcU}{\mathcal{U}}
\newcommand{\mcV}{\mathcal{V}}

\newcommand{\mcX}{\mathcal{X}}
\newcommand{\mcY}{\mathcal{Y}}

\DeclareMathOperator{\exact}{ex}

\DeclareMathOperator{\rad}{rad}
\DeclareMathOperator{\Zg}{Zg}
\DeclareMathOperator{\mSpec}{mSpec}

\DeclareMathOperator{\add}{add}

\usepackage{todonotes}

\begin{document}
\title{Classification of exact substructures using the Ziegler spectrum}
\author{Julia Sauter}
\address{Julia Sauter\\ Faculty of Mathematics \\
Bielefeld University \\
PO Box 100 131\\
D-33501 Bielefeld }
\email{jsauter@math.uni-bielefeld.de}

\begin{abstract}
Given an idempotent complete additive category, we show the there is an explicitly constructed topological space such that the lattice of exact substructures is anti-isomorphic to the lattice of closed subsets. 
In the special case that the additive category has weak cokernels, this topological space is an open subset of the Ziegler spectrum and this is a result of Kevin Schlegel.    
\end{abstract}


\maketitle

\section{The main result}

Let $\mcC$ be an idempotent complete small additive category and $\exact (\mcC )$ the lattice of all exact structures on $\mcC$. \\
Let $\mcA$ be the ind-completion of $\mcC$ and $\Zg (\mcA)^{0}$ be the set of all indecomposable pure-injective modules which are not injective in the maximal locally coherent exact structure. 
Then we construct a topology on $\Zg(\mcA)^0$ such that we have an explicit anti-isomorphism of lattices between $\exact (\mcC )$ and the closed subsets in $\Zg(\mcA)^0$.  

\begin{thm} (cf. Theorem \ref{FirstBij} and for $\mcC$ with weak cokernel, cf \cite[Thm B]{Sch}) 
Let $\mcC$ be a small idempotent complete additive category with maximal exact structure $\mcE_{max}$, then  
\[
\begin{aligned}
\{ \text{closed sets in } \Zg (\mcA)^0\} & \longleftrightarrow \exact (\mcC)\\
\mcU &\longmapsto (\mcE_{\rm max})^{\mcU} \quad \text{ and conversely } \mcE\longmapsto \mcU_{\mcE}\end{aligned}
\]
where (1) $\mcE_{max}^{\mcU}$ consists of all $\mcE_{max}$-short exact sequences $\sigma$ such that $\Hom_{\mcA} (\sigma , U)$ is exact for all $U \in \mcU$\\
and (2) $\mcU_{\mcE}$ is the set of indecomposables $\overrightarrow{\mcE}$-injectives (where $\overrightarrow{\mcE}$ is the ind-completion of the exact category $\mcE$)\\
give mutually inverse bijections. Furthermore, we have for $\mcE, \mcE' $ in $\exact (\mcC)$: 
\[
\mcE\leq \mcE' \quad \Leftrightarrow \quad \mcU_{\mcE} \supseteq \mcU_{\mcE'}
\]
\end{thm}

We remark that we do not know if there is a topology on $\Zg(\mcA)$ (if $\mcA$ does not have products). 
After proving this theorem, we recall the 'representation-finite' case (here we recall a theorem of Enomoto). Then we look at module categories of rings where the Ziegler spectrum is well-known and study the parametrizations of exact substructures explicitly (we look at a commutative discrete valuation ring, at commutative Dedekind domains and the path algebra of the Kronecker quiver). In some cases, we use exact substructures 'making a torsion functor exact' - therefore we introduce this concept in a separate subsection. \\
We also give an Appendix on Ind-completion of exact categories (summarizing results from \cite{P-locCoh}). 

\subsection*{Acknowledgement:}
This work has been supported by the Deutsche Forschungsgemeinschaft (DFG, German Research Foundation) -- Project-ID 491392403 -- TRR 358.

\section{Proof of the main result}
We refer to the Appendix for Ind-completion of small idempotent complete additive and exact categories. For an essentially small additive category $\mcC$ we call $\overrightarrow{\mcC}=:\mcA$ its Ind-completion. 
Let $\mcE$ be an exact structure and $\overrightarrow{\mcE}$ be its ind-completion and $i \colon \mcE \to \overrightarrow{\mcE}$ the Yoneda embedding. We call $X$ in $\overrightarrow{\mcE}$ is \textbf{fp}-$\overrightarrow{\mcE}$-\textbf{injective} if $\Ext^1_{\overrightarrow{\mcE}} (i(E), X)=0$ for all $E$ in $\mcE$. Let $\mcX_{\mcE}$ be the full subcategory of fp-injectives in $\overrightarrow{\mcE}$.

\begin{dfn}
An object in an additive category $M$ is called \textbf{indecomposable} if $M=N\oplus L$ implies $N=0$ or $L=0$.\\
Let $\mcA=\overrightarrow{\mcC}$ be a locally finitely presented additive category. We define the \textbf{pure exact structure} to be $\overrightarrow{\mcE_0}=\overrightarrow{\mcC_{\textbf{split}}}$ to be the ind-completion of the split exact structure. An object in $\mcA$ is called \textbf{pure-injective} if it is an injective in $\overrightarrow{\mcE_0}$ . 
Then we define the following class 
\[
\Zg (\mcA):= \{ [M]\mid M \in \mcA \text{ indecomposable and pure-injective}\}
\]
where $[M]$ denotes the isomorphism class of an object $M$. 
\end{dfn}
We will usually pick representatives for the isomorphism classes in $\Zg(\mcA)$ and treat elements in $\Zg (\mcA)$ as objects in $\mcA$. 

We recall the following definition from \cite{BeG}
\begin{dfn}
A \textbf{quasi-limit} of a (covariant) functor $H\colon I \to  \mcB$ consists of an object in $\mcB$ together with a family of morphisms \[(\lim{}^q H , \;\;\pi_i \colon \lim{}^q H \to H(i), i \in I)\] such that for every morphism $\alpha\colon  i \to j$ we have $H(\alpha )\pi_i=\pi_j$ (we say $(\pi_i)$ is an compatible family) and \\
\begin{itemize}
\item[(1)] For every compatible family $z_i \colon B \to H(i)$, $i\in I$, we have a (not necessary unique) $z\colon B \to \lim^q H$ such that $\pi_i z=z_i$ for all $i \in I$.
\item[(2)] If $h \colon \lim^qH \to \lim^qH$ is a morphism such that $\pi_i h =\pi_i$ for all $i \in I$, then $h$ is an automorphism. 
\end{itemize}
A \textbf{quasi-product} is a quasi-limit for a small  category $I$ with only identity morphisms, in this case we denote the quasi-limit as $\prod^q_{i\in I}X_i$ (where $H(i)=X_i$, $i \in I$). \\ 
We say $\mcB$ is \textbf{quasi-complete} (resp. has \textbf{quasi-products}) if every functor $H \colon I \to \mcB$ has a quasi-limit in $\mcB$ (resp. only for categories $I$ with only identity morphisms). 
\end{dfn}

\begin{lem}\cite[Thm 2]{BeG}
If $\mcA$ is locally finitely presented then it has  quasi-limits. In particular it also has quasi-products. \end{lem}
Identify $\mcA$ with ${\rm Flat}(\mcC^{op}, Ab)$ (cf. \cite[Thm 3]{CB}), this is a subcategory of $\Modu \mcC$ (all additive functors $\mcC^{op} \to (Ab)$). Then in loc. cit. quasi-limits are provided by taking the limit in $\Modu \mcC$ and then compose with its flat cover.  

\begin{dfn}
If $S$ is a set of objects in a locally finitely presented additive category, then we say $S$ is a \textbf{pure cogenerator} if every object $A$ in $\mcA$ fits into a pure exact sequence
\[
A\infl \prod{}^q_{i\in I} E_i \defl B
\]
with $E_i\in S$, $i \in I$ (and $I$ a set). 
\end{dfn}

\begin{thm} (\cite[Cor 4]{BeG})
If $\mcA$ is locally finitely presented then 
$\Zg (\mcA )$ is a set and a pure cogenerator for $\mcA$.     
\end{thm}
In particular, it is also a non-empty set.

\begin{dfn}
Let $\mcC$ be an essentially small additive category, $\mcA=\overrightarrow{\mcC}$ and $S$ a class of morphisms in $\mcC$. 
Let $\mcX(S)$ be the full subcategory of $\mcA$ of all objects $I$ with the following property: For any map $s\colon M \to M'$
in $S$ and any map $f \colon M \to I$ there exists $f'\colon M'\to I$ such that $f's=f$. Alternatively, one can descibe this as 
\[
\mcX (\mcS) =\{I\in \mcA\mid \coKer \Hom_{\mcA} (s, I)=0 \quad \forall s\in S\}
\]
A full subcategory $\mcX$ of $\overrightarrow{\mcC}$
is called \textbf{definable} if there exists a class of morphisms $S$ in $\mcC$ such that $\mcX=\mcX(S)$. 
\end{dfn}

We observe the following (straightforward):
\begin{lem} \label{firstProperties}
Let $\mcX$ be definable and $X\infl Y \defl Z$ a pure exact sequence in $\mcA$. \\
\begin{itemize}
\item[(0)] $\mcX$ is closed under direct summands. 
\item[(1)] If $X,Z$ in $\mcX$ then $Y$ in $\mcX$
\item[(2)] If $Y$ in $\mcX$ then $Z\in \mcX$
\item[(3)] If $\mcX=\mcX(S)$ and all morphisms in $S$ have weak cokernels (in $\mcC$), then $Y$ in $\mcX$ implies $X \in \mcX$. 
\end{itemize}
\end{lem}

In particular, if a definable subcategory is closed under pure subobjects then it is a pure Serre subcategory (i.e. a Serre subcategory in the pure exact structure). 

In general, we do not know when definable subcategories are closed under quasi-products, but we observe the following (which applies in several situations). 

\begin{lem} \label{qprods}
Let $\mcF$ be an exact category and assume that the underlying additive category $\mcA$ has quasi-products. Let $\mcB$ be a class of obects and \[\mcX:=\{X\in \mcA\mid \Ext^1_{\mcF} (B,X)=0 \quad \forall B \in \mcB\}.\]
Then $\mcX$ is closed under quasi-products. 
\end{lem}

\begin{proof} 
Let $X_i$, $i \in I$ be a family of objects in $\mcX$ and $X=\prod{}^q_{i\in I}X_i$, we denote the canonical morphisms by $\pi_i \colon X \to X_i$, $i \in I$. Assume we have an $\mcF$-exact sequence 
\[
\begin{tikzcd}
X \arrow[r, "g", tail] & Y \arrow[r, two heads] & B
\end{tikzcd}
\]
with $B$ in $\mcB$. By assumption, when we form the pushout along $\pi_i$, we get a split exact sequence, so we find $g_i \colon X_i \to Y_i $ and $r_i \colon Y_i\to X_i$ with $r_ig_i={\rm id}_{X_i}$ and a commuting diagram 
\[
\xymatrix{ X\ar[r]^g \ar[d]_{\pi_i}& Y\ar[d]_{y_i} \\
X_i \ar[r]_{g_i}& Y_i
}
\]
Now, the universal property of the quasi-product applied to the family $r_i y_i \colon Y \to X_i$, $i \in I$ implies the existence of a morphism $h\colon Y \to X$. Now, we look at the endomorphism $hg\colon X \to X$. Then we have $\pi_i hg = \pi_i$ for all $i \in I$ and using the minimality property of the quasi-product we conclude that $hg$ is an isomorphism and therefore $g$ a split monomorphism. 
\end{proof}

Here are the examples we would like to apply it to:
\begin{exa}\label{Exqprods}
Let $\overrightarrow{\mcE}$ be a locally coherent exact category. \\
The subcategory $\mcI(\overrightarrow{\mcE})$ of injectives in $\overrightarrow{\mcE}$ is closed under quasi-products. \\
In particular, the subcategory of pure-injectives in a locally finitely presented additive category is always closed under quasi-products. \\
Also, the subcategoy $\mcX_{\mcE}$ is closed under quasi-products.
\end{exa}

We recall a concept of Herzog (cf. \cite{Her})
\begin{dfn} 
Let $\mcD=\{X \mid \Ext^1_{\Modu \mcC} (F,X)=0 \; \forall F \in {\rm Flat}(\mcC^{op}, Ab)\}\subseteq \Modu \mcC$ (this category is called the category of \textbf{cotorsion modules}). Let $M\in \Modu \mcC$, then a $\mcD$-\textbf{preenvelope} is given by a monomorphism $i$ 
\[
M \xrightarrow{i} {\rm CE}(M) \defl Z=\coKer (i)
\]
such that ${\rm CE}(M)$ in $\mcD$ and $Z$ flat (i.e. in ${\rm Flat}(\mcC^{op}, Ab)$). It is called a $\mcD$-\textbf{envelope} if additionally: Whenever $t \colon {\rm CE}(M)\to {\rm CE}(M)$ fulfills $ti=i$, then $t$ is an automorphism. \\
If $M$ is flat, then also ${\rm CE}(M)$ and therefore it is pure-injective. In this case, ${\rm CE}(M)$ is called the \textbf{pure-injective envelope} of $M$ and $M \infl {\rm CE} (M) \defl Z$ is a pure exact sequence in $\mcA$.  
\end{dfn}

\begin{thm} (\cite[Thm 6]{Her}) Every object in a locally finitely presented additive category admits a pure-injective envelope. 
\end{thm}

\begin{lem} (generalization of \cite[Thm 3, Thm 4]{BeG})\label{pureCogen}
If $\mcX\subseteq \mcA$ is a pure Serre subcategory which is also closed under taking pure-injective envelopes and quasi-products, then 
\[
\Zg(\mcA) \cap \mcX :=\{ M\in \Zg(\mcA)\mid M\text{ in }\mcX\} \quad \subseteq \Zg(\mcA)
\]
is a pure cogenerator for $\mcX$, i.e. every $X$ in $\mcX$ there exists a pure exact sequence 
\[
X\infl \prod{}^q_{i\in I} E_i \defl Y
\]
with $E_i\in \Zg(\mcA) \cap \mcX$, $i \in I$ (and $I$ a set). 
\end{lem}

\begin{proof}
We identify $\mcA={\rm Flat} (\mcC^{op}, Ab)$ (cf. \cite[Thm 3]{CB}). 
We observe that the (construction in the) proof of \cite[Thm 3]{BeG} generalizes to the following: If $F$ is in $\mcX$, $x\in F(C)$ a non-zero element. Then there exists a $E_x\in \mcX\cap \Zg(\mcA)$ and a natural transformation $\tau_x\colon F\to E_x$ such that $\tau_x(x)\neq 0$ in $E_x(C)$.\\ 
Then using the proof of \cite[Thm 4]{BeG} the claim follows.       
\end{proof}

\begin{dfn}
We call a subcategory $\mcX$ of $\mcA$ \textbf{strongly definable} if it is definable and closed under pure subobjects, quasi-products and pure-injective envelopes. \\
A subset $\mcU \subseteq \Zg (\mcA)$ is called \textbf{Ziegler-closed} if there exists a \emph{strongly} definable subcategory $\mcX$ such that $\mcU= \Zg (\mcA)\cap \mcX$.     
\end{dfn}

\begin{cor} \label{bijStronglyDef}
The map $\mcX\mapsto \mcX\cap \Zg(\mcA)$ gives a bijection strongly definable subcategories and Ziegler-closed subsets.     
\end{cor}

This corollary is obvious: Given $\mcX$ strongly definable, by Lemma \ref{pureCogen}, the set $\Zg(\mcA)\cap \mcX$ is a pure cogenerator of $\mcX$. As $\mcX$ is a pure Serre subcategory, this implies that the pure cogenerator determines all $\mcX$.

\begin{thm} \label{FirstBij}
Let $\mcC$ be a small idempotent complete additive category and $\mcE_{\max}$ its maximal exact structure. For every small exact category $\mcE$ let us denote by $\mcU_{\mcE}={\rm indec} (\mcI(\overrightarrow{\mcE}))$ the set of all indecomposable objects in the category of all $\overrightarrow{\mcE}$-injectives and we write 
$\mcU_{max}:=\mcU_{\mcE_{max}}$. \\
Then the following maps give mutually inverse bijections  
\[
\begin{aligned}
\{\mcU \text{ Ziegler-closed }, \mcU_{max} \subseteq \mcU\} & \longleftrightarrow \exact (\mcC)\\
\mcU &\longmapsto (\mcE_{\rm max})^{\mcU} \quad \text{with inverse } \; \mcE\longmapsto \mcU_{\mcE}\end{aligned}
\]
where $\mcE_{max}^{\mcU}$ consists of all $\mcE_{max}$-short exact sequences $\sigma$ such that $\Hom_{\mcA} (\sigma , U)$ is exact for all $U \in \mcU$. 
\end{thm}

We start with the well-definedness of these assignments, we first need: 
\begin{pro}
The category $\mcX_{\mcE}$ of fp-$\overrightarrow{\mcE}$-injectives is strongly definable.     
\end{pro}
We show it is definable and closed under pure subobjects: 
\begin{lem} \label{fp-injsDefinable}Let 
Let $\mcC$ be essentially small, idempotent complete  and $\mcE$ an exact structure on it. Let $\overrightarrow{\mcE}$ be its Ind-completion. 
Let $\mcX_{\mcE}$ be the full subcategory fp-$\overrightarrow{\mcE}$-injectives, then this is a strongly definable subcategory since it can be written as 
\[
\mcX_{\mcE}=\mcX({\rm Infl}_{\mcE})
\]
where ${\rm Infl}_{\mcE}$ denotes the $\mcE$-inflations.    
\end{lem}

We remark the following straightforward corollary of this Lemma: 
\begin{cor}
The subcategory $\mcX_{\mcE}$ is inflation-closed in $\overrightarrow{\mcE}$. Since $\overrightarrow{\mcE}$ has enough injectives and as $\mcX_{\mcE}$ contains the injectives in $\overrightarrow{\mcE}$, is closed under direct summands and extensions, it follows:\\
$\mcX_{\mcE}$ is a coresolving subcategory in $\overrightarrow{\mcE}$.
\end{cor}

\begin{proof} (of Lemma \ref{fp-injsDefinable}) 
By definition we have that $\mcX ({\rm Infl}_{\mcE})$ is the full subcategory of all objects $I$ in $\mcA$ such that $\Hom (\sigma , I)$ is exact for all $\mcE$-short exact sequences. These objects clearly contain the fp-$\overrightarrow{\mcE}$-injectives. We show they are equal. Let $X\in \mcX({\rm Infl}_{\mcE})$, $E\in \mcE$ and given $\sigma$ a short exact sequence 
$X\infl J \defl E$ in $\overrightarrow{\mcE}$. By \cite[Lem 1.5]{P-locCoh} (take as morphism the identity in $E$), there exists a commutative diagram 
\[
\begin{tikzcd}
X' \arrow[r, tail] \arrow[d] & J' \arrow[r, two heads] \arrow[d] & E \arrow[d, "\rm id"] \\
X \arrow[r, tail]            & J \arrow[r, two heads]            & E                    
\end{tikzcd}
\]
with the upper row a short exact sequence $\sigma'$ in $\mcE={\rm fp}(\overrightarrow{\mcE})$. By \cite[Prop 2.12]{Bue} we have that $\sigma$ is the pushout of $\sigma'$ along $X'\to X$. By assumption, we have $\Hom (\sigma', X)$ is exact. As the condition $\Hom (-,X)$ being exact is closed under forming push-out, this imples $\Hom (\sigma , X)$ is exact. Then looking at this exact sequence, we see that ${\rm id} \in \Hom (X,X)$ has a preimage in $\Hom (J,X)$ but this translates into $X\infl J$ is a split monomorphism. We conclude that $\Ext^1_{\overrightarrow{\mcE}}(E,X)=0$.   
\end{proof}

We we can conclude the previous Proposition with: 
\begin{lem} \label{PI-envelopes}
The subcategory $\mcX_{\mcE}$ is closed under pure-injective envelopes, in particular it is strongly definable.      
\end{lem}

The following observation is important to understand the lemma.
\begin{rem} \label{key}
By the definition of $\mcX_{\mcE}$ as fp-$\overrightarrow{\mcE}$-injectives: Every $\overrightarrow{\mcE}$-exact sequence $X \infl Y \defl Z$ with $X$ in $\mcX_{\mcE}$ is pure exact.    
\end{rem}

\begin{proof}
Let $X$ be an object in $\mcX_{\mcE}$. We look at the inclusions $\mcX_{\mcE} \subseteq \mcA\cong {\rm Flat} (\mcC^{op}, {\rm Ab})\subseteq \Modu \mcC$. 
The pure-injective envelope is the $\mcD$-envelope in $\Modu-\mcC$ with $\mcD:=\{M\in \Modu \mcC\mid \Ext^1_{\Modu \mcC} (F,M)=0 \; \forall F \in {\rm Flat} (\mcC^{op}, {\rm Ab})\}$. The category $\mcD$ is called the category of cotorsion modules, and $X\infl {\rm CE}(X) \defl Z$ be the defining short exact sequence in $\Modu \mcC$ for the $\mcD$-envelope. \\
Now, we also know that we have an $\overrightarrow{\mcE}$-exact sequence $X\infl E(X) \defl Y$ with $E(X)$ in $\mcI(\mcE)$, by remark \ref{key} this is a pure exact sequence and $E(X)$ is pure-injective and lies in $\mcX_{\mcE}$.
We look at the inclusion $\overrightarrow{\mcE}={\rm Flat} (\mcC^{op}, \rm Ab)\subseteq \Modu \mcC$. 
As $E(X)$ is a cotorsion module we get a commutative diagram (by the property of the $\mcD$-preenvelope)
\[
\xymatrix{
X \ar[r]^{j}\ar[dr]^{i}& {\rm CE}(X)\ar[d]^{d} \\
& E(X)
}
\]
Now, the exact sequence $X \infl E(X) \defl Y$ is pure exact, ${\rm CE}(X)$ is pure-injective, so when we apply $\Hom (-, {\rm CE}(X))$ to leads it gives a surjective map $\Hom (E(X), {\rm CE}(X)) \to \Hom (X, {\rm CE}(X))$. In particular there exists a morphism $g\colon E(X) \to {\rm CE}(X)$ such that $gi=j$. Now, we find a commuting diagram 
\[
\xymatrix{
X\ar[r]^j\ar[dr]_j & {\rm CE}(X)\ar[d]^{gd} \\
& {\rm CE}(X)},\]
and so by the property of a $\mcD$-envelope, $gd$ is an isomorphism and then $d$ is a split monomorphis. As $\mcX_{\mcE}$ is closed under direct summands, it follows that ${\rm CE}(X)$ lies in $\mcX_{\mcE}$. 
\end{proof}

\begin{lem} \label{basicDef} 
Let $\mcC$ be essentially small, idempotent complete  and $\mcE$ an exact structure on it. Let $\overrightarrow{\mcE}$ be its Ind-completion. 
Let $\mcU_{\mcE}={\rm indec } (\mcI(\overrightarrow{\mcE}))$ be the set of indecomposable injectives in $\overrightarrow{\mcE}$. Then 
\[
\mcU_{\mcE}=\Zg (\mcA) \cap \mcX_{\mcE}
\]
and therefore $\mcU_{\mcE}$ is Ziegler-closed. 
\end{lem}

\begin{proof}
Clearly $\mcU_{\mcE} \subseteq \Zg (\mcA) \cap \mcX_{\mcE}$, we need to see the other inclusion.\\ 
Let $U\in \Zg (\mcA) \cap \mcX_{\mcE}$. We consider $\sigma$ a short exact sequence $U \infl Y\defl X$ in $\overrightarrow{\mcE}$, as $U$ is pure-injective, we only need to see that $\sigma$ is pure exact (for $\sigma$ to be split exact). We write $\sigma=\colim \sigma_i$ as a filtered colimit of $\mcE^{\mcU}$-short exact sequences $U_i \infl Y_i \defl X_i$. Now, we factorize the canonical morphisms $\sigma_i \to \sigma$, $i \in I$ of short exact sequences following \cite[Prop. 3.1]{Bue}
\[
\begin{tikzcd}
\sigma_i  \arrow[d] & U_i \arrow[r, tail] \arrow[d] & Y_i \arrow[r, two heads] \arrow[d] & X_i \arrow[d,-,double equal sign distance,double] \\
\eta_i \arrow[d]    & U \arrow[r, tail] \arrow[d,-,double equal sign distance,double]   & Z_i \arrow[r, two heads] \arrow[d] & X_i \arrow[d] \\
\sigma              & U \arrow[r, tail]             & Y \arrow[r, two heads]             & X            
\end{tikzcd}
\]
This means $\eta_i$ is the push-out of $\sigma_i$ along the canonical morphism $U_i \to U$. As $\Hom_{\mcA} (\sigma_i, U)$ is exact, it follows that $\eta_i$ is split exact. Now, it is a straight forward observation to see that we have $\colim_I \eta_i =\sigma$. As $\eta_i$ are split exact they are also pure exact sequences. Now, filtered colimits of pure exact sequences are again pure exact as the pure exact structure is a locally coherent exact structure (cf. Appendix \ref{UElocCoh}). In particular $\sigma$ is pure exact. This implies $U$ in $\mcU_{\mcE}$. 
\end{proof}

We need the following easy lemma for the proof.
\begin{lem}\label{fp-inj}
Let $\mcC$ be an idempotent complete essentially small category and $\mcA=\overrightarrow{\mcC}$ its Ind-completion. Assume we have an exact structure $\mcE$ on $\mcC$ and $\mcU$ some set of pure-injective objects in $\mcA$. We denote $\mcE^{\mcU}$ the exact substructure of $\mcE$ consisting of $\mcE$-exact sequences $\sigma$ such that $\Hom_{\mcA} (\sigma , U)$ exact for all $U \in \mcU$. 
Then all objects in $\mcU$ are ${\rm fp} (\overrightarrow{\mcE^{\mcU})}$-injectives. 
\end{lem}

\begin{proof} (of Lemma \ref{fp-inj}) The proof follows by the same argument as in the proof of Lemma \ref{basicDef}. Let $U$ be in $\mcU$ and $X$ be an $\mcC$ and we take a $\overrightarrow{\mcE^{\mcU}}$-short exact sequence 
\[
\sigma \colon U \infl Y \defl X
\]
We need to see it splits. We write $\sigma=\colim \sigma_i$ as a filtered colimit of $\mcE^{\mcU}$-short exact sequences $U_i \infl Y_i \defl X_i$. Now, we factorize the canonical morphisms $\sigma_i \to \sigma$, $i \in I$ of short exact sequences following \cite[Prop. 3.1]{Bue}
\[
\begin{tikzcd}
\sigma_i  \arrow[d] & U_i \arrow[r, tail] \arrow[d] & Y_i \arrow[r, two heads] \arrow[d] & X_i \arrow[d,-,double equal sign distance,double] \\
\eta_i \arrow[d]    & U \arrow[r, tail] \arrow[d,-,double equal sign distance,double]   & Z_i \arrow[r, two heads] \arrow[d] & X_i \arrow[d] \\
\sigma              & U \arrow[r, tail]             & Y \arrow[r, two heads]             & X            
\end{tikzcd}
\]
This means $\eta_i$ is the push-out of $\sigma_i$ along the canonical morphism $U_i \to U$. As $\Hom_{\mcA} (\sigma_i, U)$ is exact, it follows that $\eta_i$ is split exact. Now, it is a straight forward observation to see that we have $\colim_I \eta_i =\sigma$. As $\eta_i$ are split exact they are also pure exact sequences. Filtered colimits of pure exact sequences are again pure exact as the pure exact structure is a locally coherent exact structure (cf. Appendix \ref{UElocCoh}). In particular, $\sigma$ is pure exact and $U$ is pure-injective, therefore it splits.   
\end{proof}

Let us come back to: 
\begin{proof} (of Thm. \ref{FirstBij})
Let $\mcE$ be an exact structure on $\mcC$ and we set $\mcU:=\mcU_{\mcE}$. As $\mcE$ is fully exact in $\overrightarrow{\mcE}$ we have that $\mcE\leq \mcF:= \mcE_{max}^{\mcU}$ is an exact substructure. To see that they are equal, it is enough to see that $\mcU=\mcU_{\mcF}$. As $\mcE\leq \mcF$ we have that $\overrightarrow{\mcE}\leq \overrightarrow{\mcF}$ this implies that the subcategory of injectives fulfill $\mcI(\overrightarrow{\mcE})\supseteq \mcI (\overrightarrow{\mcF})$ and therefore $\mcU \supseteq \mcU_{\mcF}$.  Now, for the other inclusion we conclude from Lemma \ref{fp-inj} that $\mcU\subseteq \mcX_{\mcF}$. This implies $\mcU \subseteq \mcX_{\mcF} \cap \Zg (\mcA) = \mcU_{\mcF}$ by Lemma \ref{basicDef}. To finish the proof it is enough to see that Ziegler-closed $\mcU$ containing $\mcU_{max}$ are always of the form $\mcU_{\mcE}$ for an exact structure $\mcE$ on the additive category $\mcC$. We discuss this separately in the next propostion. 
\end{proof}
Recall first the following result of Enomoto: 
\begin{thm} (\cite[Thm 2.7]{E})\label{EnoBij}
Let $\mcC$ be idempotent complete small additive category with maximal exact structure $\mcE_{max}$. Then $\eff(\mcE_{\max})$ is a Serre subcategory in $\modu_1\mcC$. 
The assignment $\mcE\mapsto {\rm eff}(\mcE)$ gives a bijective map from (1) to (2) with 
\begin{itemize}
\item[(1)] exact structures on $\mcC$
\item[(2)] Serre subcategories of ${\rm eff}(\mcE_{max})$
\end{itemize}
We consider the set of Serre subcategories as a poset with respect to inclusion. Then this bijection is an isomorphism of posets. 
\end{thm}
This version is not the given one in loc.cit. but combine it with: Enomoto constructed in \cite[Lem 2.5]{E} a duality which restricts to a duality $E\colon \eff(\mcE_{max})\to \eff (\mcE_{max}^{op})$. As both are abelian categories, it maps Serre subcategories to Serre subcategories. This implies that the Serre subcategories in $\eff(\mcE_{max})$ are those given in loc. cit.

Let $s\colon M \to M'$ be a morphism in $\mcC$.
Then the functor $F=F_s=\coKer \Hom_{\mcC} (s,-)\colon \mcC \to (Ab)$ extends to a functor $\overline{F}=\overline{F}_s=\coKer \Hom_{\mcA} (s,-)\colon \mcA \to (Ab)$. We have for a filtered colimit $A=\colim_{i\in I} X_i$ with $X_i\in \mcC$ that $\Hom_{\mcA} (C, A)=\colim_{i\in I} \Hom_{\mcC} (C, X_i)$ for all $C$ in $\mcC$ since $\mcC={\rm fp}(\mcA)$. \\
Also since $(Ab)= \Modu \mathbb{Z}$ is a locally coherent abelian category, taking filtered colimits of short exact sequence is exact and we can conclude that 
\[
\overline{F}(A) = \coKer \Hom_{\mcA}(s,A)= \colim_{i \in I} \coKer \Hom_{\mcC} (s, X_i) = \colim_{i \in I} F(X_i)
\]
Then a definable subcategory is $\mcX(S)=\bigcap_{s\in S} \Ker \overline{F}_s$.\\
Observe that the functor $F\mapsto \overline{F}$ has the following property, let $A=\colim_{i\in I}X_i$ be a filtered colimit with $X_i \in \mcC$, then for every short exact sequence $F\infl G \defl H$ in $\modu_1(\mcC^{op})$ we have short exact sequence in $(Ab)$
\[
F(X_i) \infl G(X_i) \defl H(X_i) , \quad i \in I
\]
and as filtered colimits are exact in $(Ab)$ we find a short exact sequence of abelian groups  
\[
\overline{F} (A) \infl \overline{G}(A) \defl \overline{H} (A) 
\]
\\
Then, clearly given $\mcX\subseteq \mcA$, the subcategory  
\[
L_{\mcX} := \{F\in \modu_1 (\mcC^{op}) \mid \mcX\subseteq \Ker \overline{F} \}
\]
defines a Serre subcategory in $\modu_1 (\mcC^{op})$.
Conversely, every Serre subcategory $L$ in $\modu_1 (\mcC^{op})$ defines a definable subcategory via 
$\mcX^L:=\bigcap_{F \in L} \Ker \overline{F}$.

The following bijection is a trivial generalization from the previously known bijection: 
\begin{lem} \label{SerreBij}
The assignments $\mcX\mapsto L_{\mcX}$ and $L \mapsto \mcX^L$ are mutually inverse bijections between 
\begin{itemize}
\item[(1)] definable subcategories $\mcX$ in $\mcA$
\item[(2)] Serre subcategories $L$ of $\modu_1(\mcC^{op})$.
\end{itemize}
This bijection fulfills for two definable subcategories $\mcX, \mcX'$:\\
$\mcX\subseteq \mcX'$ iff and only if $L_{\mcX} \supseteq L_{\mcX'}$. 
\end{lem}

\begin{proof}
If $\mcX$ is definable, then it is straight-forward to see
\[
\mcX= \bigcap_{F \in L_{\mcX}} \Ker \overline{F}
\]
If $L$ is a Serre subcategory, we want to see $L=L_{\mcX^{L}}$. Clearly $L\subseteq  L_{\mcX}$. For the other inclusion: Assume $F\notin L_{\mcX^L}$, then $\mcX^L$ is not a subcategory of $\Ker \overline{F}$. This means there exists an object $A$ in $\mcA$ such that $\overline{F} (A)\neq 0$ and $\overline{H} (A)=0$ for all $H$ in $L$. In particular $\overline{F}\neq \overline{H}$ for all $H$ in $L$ and therefore $F\notin L$. \\
The last statement is very easy to see. 
\end{proof}

\begin{exa}
$\mcX_{\mcE}$ corresponds to the Serre subcategory ${\rm eff}(\mcE^{op})$. 
\end{exa}

\begin{cor}
Let $\mcU=\mcX^L\cap \Zg(\mcA)\subseteq \Zg (\mcA )$ be Ziegler-closed for some Serre subcategory $L$ such that $\mcX^L$ is strongly definable with $\mcU_{max}\subseteq \mcU$.\\
Then $\mcX^{\eff(\mcE^{op}_{max})} \subseteq \mcX^L$ by Corollary \ref{bijStronglyDef} and then 
$L\subseteq {\rm eff} (\mcE^{op}_{max})$ by Lemma \ref{SerreBij}\\
In particular, there exists an exact structure $\mcE$ such that $L={\rm eff }(\mcE^{op})$ by Thm \ref{EnoBij} and so $\mcU=\mcU_{\mcE}$. 
\end{cor}

\begin{rem} \label{poset}
We observe the following properties of the bijection:
Assume $\mcE=\mcE^{\mcU}$, $\mcF=\mcE^{\mcV}$ are two exact structures on $\mcC$ corresponding to two Ziegler-closed subsets $\mcU_{\max} \subseteq \mcU, \mcV \subseteq \Zg (\mcA)$. Then: \\
\[\mcE\leq \mcF \quad \Leftrightarrow \quad \mcU \supseteq \mcV \]
\end{rem}

Let us note the following result:
\begin{pro} \label{keyForTop}
If $\mcE$ and $\mcF$ are two exact structures on a small idempotent complete additive category $\mcC$, recall that $\mcE \wedge \mcF$ denotes its meet in the lattice of exact structures. Then we have 
\[
\mcU_{\mcE\wedge \mcF} =\mcU_{\mcE}\cup \mcU_{\mcF}
\]    
\end{pro}

We need the following two lemmata for the proof: 
\begin{lem} \label{End}(\cite[Cor 1]{BeG}) Let $\mcA$ be locally finitely presented. 
Let $M$ be an indecomposable pure-injective in $\mcA$, then 
$S=\End_{\mcA} (M)$ is a local ring. 
\end{lem}
We remark the following about the proof: 
\begin{rem} 
This Lemma is \cite{BeG}, Cor 1, but as W. Crawley-Boevey observed: It is a corollary of Proposition 3 in loc. cit. which contains a wrong statement (not $S$ is regular but $S/J(S)$ is regular (they misquoted Asensio-Herzog). Nevertheless, \cite{BeG} show in the proof of Prop. 3 that $S$ is an indecomposable right cotorsion ring and this implies that $S$ is local.  
\end{rem}

\begin{lem} \label{meets}
Let $\mcE$ and $\mcF$ be two exact structures on a small additive category $\mcC$. Then:   
\[\overrightarrow{\mcE}\wedge \overrightarrow{\mcF} =\overrightarrow{\mcE\wedge\mcF}\]
\end{lem}
\begin{proof} (of Lemma \ref{meets}) Clearly $\overrightarrow{\mcE\wedge \mcF}\leq \overrightarrow{\mcE}\wedge \overrightarrow{\mcF}$. 
We use \cite[Lemma 1.5]{P-locCoh} to show that every 
$\overrightarrow{\mcE}\wedge \overrightarrow{\mcF}$-deflation $d\colon E\defl X$ is also a $\overrightarrow{\mcE\wedge \mcF}$-deflation. Let $f\colon S\to X$ be a morphism with $S\in \mcC$. By loc. cit. there exists an $\mcE$-deflation $u \colon U \to S$ and an $\mcF$-deflation $v\colon V \to S$ and morphisms $e\colon U \to E, g\colon V \to E$ such that 
$fu=de$ and $fv=dg$. Now, $\left(\begin{smallmatrix} v\\u \end{smallmatrix} \right)\colon V\oplus U \to S$ is an $\mcE\wedge \mcF$-deflation (use \cite[Prop. 2.12]{Bue}) and we have a commuting diagram 
\[
\begin{tikzcd}
V\oplus U \arrow[r, "\left(\begin{smallmatrix} v \\u \end{smallmatrix}\right)"] \arrow[d, "\left(\begin{smallmatrix}g\\e\end{smallmatrix}\right)"'] & S \arrow[d, "f"] \\
E \arrow[r, "d"']                                                                                                                                  & X               
\end{tikzcd}
\]
Then by \cite[Lem. 1.5]{P-locCoh} we conclude that $d$ is an $\overrightarrow{\mcE\wedge \mcF}$-deflation.
\end{proof}

\begin{proof} (of Prop. \ref{keyForTop})
By remark \ref{poset}, we conclude 
\[\mcU_{\mcE\wedge \mcF} \supseteq \mcU_{\mcE}\cup \mcU_{\mcF}, \]
so we need to see the other inclusion. For $M\in \Zg(\mcA)$ we show: $M \notin \mcU_{\mcE}\cup \mcU_{\mcF}$ implies $M\notin \mcU_{\mcE\wedge \mcF}$.\\
So assume that there exists non-split $\overrightarrow{\mcE}$-exact sequence $\sigma\colon M \infl L \defl X$ and a non-split $\overrightarrow{\mcF}$-exact sequence $\tau\colon M \infl N \defl Y$, we form the push-out alongside the inflations 
\[
\begin{tikzcd}
M \arrow[r, "a", tail] \arrow[d, "b"', tail]          & L \arrow[d, tail] \arrow[r, two heads]      & X \arrow[d, "\rm{id}"] \\
N \arrow[r, tail] \arrow[d, two heads]                & E \arrow[r, two heads] \arrow[d, two heads] & X                      \\
Y \arrow[r, phantom, shift left] \arrow[r, "\rm{id}"] & Y                                           &                       
\end{tikzcd}
\]
By \cite[Prop. 2.12]{Bue} we get an induced $\overrightarrow{\mcE}\wedge\overrightarrow{\mcF}$-exact sequence $\eta$
\[
\begin{tikzcd}
M \arrow[r, "\left(\begin{smallmatrix}b \\a\end{smallmatrix}\right)", tail] & L\oplus N \arrow[r, two heads] & E
\end{tikzcd}
\]
By Lemma \ref{meets}, we have $\overrightarrow{\mcE}\wedge\overrightarrow{\mcF}=\overrightarrow{\mcE\wedge \mcF}$. We claim that $\eta$ is non-split. Assume that it is split, i.e. we find $(c, d)\colon L \oplus N \to M$ such that 
$1_M=ca+db$. By Lemma \ref{End}, we have that $S=\End_{\mcA} (M)$ is local (i.e. the non-invertible elements form an ideal). As $a$ is not a split monomorphism, the endomorphism $ca\in S$ is non-invertible. Anologously, $db\in S$ is non-invertible. Since $S$ is local, $ac+db$ non-invertible which is a contradiction. Therefore $\eta $ is non-split. 
\end{proof}

\begin{rem}
The same proof as in Lemma \ref{meets} also shows that $\overrightarrow{\bigwedge \mcE_i}=\bigwedge \overrightarrow{\mcE}_i$, but observe that the proof of Prop. \ref{keyForTop} does not generalize to infinite meets.  
\end{rem}

\begin{dfn}
A \textbf{spatial coframe} is a lattice which is isomorphic to the set of closed subsets in a topological space. The opposite lattice of a spatial coframe is called a \textbf{spatial frame} (and is more often considered).      
\end{dfn}

Then we have the following corollary of the corollary. 
\begin{cor}
We define for a subset $\mcU$ of $\Zg\setminus \mcU_{\max}$ another subset  
\[
\overline{\mcU}:=\mcU_{\mcE_{max}^{\mcU}}\setminus \mcU_{max}
\]
This operator defines a closure operator and the Ziegler-closed sets in $\Zg\setminus \mcU_{max}$ form the closed sets in a topology. \\
In particular, for every idempotent complete small additive category $\mcC$, we find a lattice isomorphism 
\[
\exact (\mcC) \to (\{ \text{ closed subsets of } \Zg(\mcA) \setminus \mcU_{max}\})^{op}
\]
where $()^{op}$ is the opposite lattice structure. 
In particular, ${\rm ex} (\mcC)$ is a spatial frame.    
\end{cor}

\begin{proof}
We need to check the following four axioms 
\begin{itemize}
\item[(1)] $\overline{\emptyset}=\emptyset$
\item[(2)] $\mcU \subseteq \overline{\mcU}$
\item[(3)] $\overline{\overline{\mcU}}=\overline{\mcU}$
\item[(4)] $\overline{\mcU \cup \mcV}=\overline{\mcU}\cup \overline{\mcV}$
\end{itemize}
Property (1) is clear, (2) follows from Lemma \ref{fp-inj}, (3) follows from the bijection in Theorem \ref{FirstBij}. 
Let us look at (4). Let $\mcE=\mcE^{\mcU}, \mcF=\mcE^{\mcV}$. By Prop. \ref{keyForTop} we have 
$\mcU_{\mcE\wedge \mcF} =\mcU_{\mcE} \cup \mcU_{\mcF}$. 
When we intersect with the complement of $\mcU_{max}$ this gives (4). 
\end{proof}

Let us also note the following corollary.
\begin{cor} Let $\mcC$ be an idempotent complete essentially small category and $\mcA=\overrightarrow{\mcC}$ its Ind-completion. Assume we have an exact structure $\mcE$ on $\mcC$ and $\mcU$ some set of pure-injective objects in $\mcA$. 
Then \[
\mcE^{\mcU} = \mcE^{\overline{\mcU\setminus \mcU_{max}}}. 
\] 
\end{cor}

\begin{proof}
By definition $\mcE^{\overline{\mcU}}$ is an exact substructure of $\mcE^{\mcU}$. This implies that $\mcU \subseteq \mcX_{\mcE^{\mcU}} \subseteq \mcX_{\mcE^{\overline{\mcU}}}$. This implies $\mcU\subseteq \mcU_{\mcE^{\mcU}}\subseteq \mcU_{\mcE^{\overline{\mcU}}}=\overline{\mcU}$ but as $\mcU_{\mcE^{\mcU}}$ is Ziegler-closed, it has to be equal to $\overline{\mcU}$. Since we have a bijection it follows that $\mcE^{\overline{\mcU}}=\mcE^{\mcU}$. 
\end{proof}

We conclude this section with the two (obviously) open  questions: 
\begin{itemize}
\item[(1)] Can we describe the Serre subcategories of $\modu_1(\mcC^{op})$ which corresponds to strongly definable subcategories (in the bijection of Lemma \ref{SerreBij})?
\item[(2)] Do all Ziegler-closed sets in $\Zg (\mcA)$ form the closed sets of a topology?
\end{itemize}

\section{The representation-finite case - Enomoto's result}

\begin{dfn} 
For a ring $\Gamma$, we denote by $\rm{proj}(\Gamma)$ the category of finitely generated projective left $\Gamma$-modules. 
We say an idempotent complete additive category $\mcA$ is \textbf{representation-finite} if it is equivalent to $\rm{proj}(\Gamma)$ for some ring $\Gamma$. 
\end{dfn}

By definition, $\mcA=\rm{proj} (\Gamma)$ is Krull-Schmidt if and only if $\Gamma$ is semi-perfect. 

We recall Enomoto's results from \cite[section 3]{E}.

\begin{lem} Let $\mcA=\rm{proj} (\Gamma )$ a Krull-Schmidt category and $\mcE$ an exact structure on $\mcE$. Then there exists an idempotent $e\in \Gamma$ such that $\mcP(\mcE) = \add (\Gamma e)$. Then $\mcE$ has enough projectives if and only if $\Gamma / \Gamma e \Gamma$ is a finite length left $\Gamma$-module. 
\end{lem}

Assume additionally that we have a commutative artinian ring $R$ and $\Gamma$ is a finitely generated $R$-algebra with $R \subset Z(\Gamma)$. This is saying that $\mcA=\rm{proj} (\Gamma)$ is Hom-noetherian $R$-linear. Then every exact structure on $\mcA= \rm{proj} (\Gamma )$ has enough projectives and enough injectives. 

\begin{dfn}
Let $\mcE$ be a an exact category and $\mcM$ a full subcategory. We say $\mcM$ is a \textbf{generator} if $\mcM$ is additively closed and for every $X$ in $\mcE$ there exists a short exact sequence $Y\to M \to X$ with $M$ in $\mcM$. A \textbf{cogenerator} is a generator in $\mcE^{op}$. \\
If $\mcA$ is an additive category and $\mcM$ a full subcategory, then we call $\mcM$ a generator (resp. cogenerator) if it is one in the maximal exact structure on $\mcA$. 
\end{dfn} 

\begin{thm} \label{Enomoto-ThmOne}(Enomoto's Theorem)
Let $R$ be a commutative artinian ring. 
Let $\mcA$ idempotent complete, representation-finite, Hom-noetherian $R$-linear. Let $\mcP=\add (\Gamma e)$ be the projectives in the maximal exact structure on $\mcA$. Then generators in $\mcA$ are given by the Boolean lattice of all additively closed subcategories containing $\mcP$, we denote it by $\rm{Generators} (\mcA)$.
Then \[
\begin{aligned}
\exact (\mcA) &\to \rm{Generators} (\mcA), \\
\mcE &\mapsto \mcP (\mcE)
\end{aligned}
\] is an isomorphism of lattices. 
\end{thm}

\begin{exa}
We look at the quiver $1\to 2 \to 3$ and at its Auslander-Reiten quiver 
\[
\begin{tikzcd}
\  & & & P_1 \arrow[rd] & &     \\
   & & P_2 \arrow[ru] \arrow[rd] \arrow[rr, dashed, dash] &  & I_2 \arrow[rd] &     \\
   & P_3 \arrow[ru] \arrow[rr, dashed, dash] &  & S_2 \arrow[ru] \arrow[rr, dashed, dash] &   & I_1
\end{tikzcd}
\]
To see the generators, fix the projectives $P_1,P_2,P_3$ and add any subset of $\{S_2,I_2,I_1\}$. So this is just the power set of this set with three elements. More interesting is to observe that we have seven hereditary exact substructures and one exact substructure of global dimension $2$, corresponding to the generator $P_1\oplus P_2\oplus P_3\oplus I_2$. 
\end{exa}

More generally for type $\mathbb{A}_n$-equioriented quivers the maximal global dimension is $n-1$. 
\begin{rem}
Exact substructures of finite-dimensional Dynkin quiver representations are always of finite global dimension. This follows because for every module $M$ over a Dynkin quiver the Gabriel of its endomorphism ring is acyclic (since the Auslander Reiten quiver of the Dynkin quiver is acyclic). This implies that the global dimension of such an endomorphism ring is finite. Then by result \cite{ASoII}, 
this implies that all exact substructures of finite global dimension. 
\end{rem}

If we are looking at representation-finite finite-dimensional algebras of global dimension $2$, then we can already find examples of exact substructures of infinite global dimension. 
\begin{exa}
Let $\Gamma$ be the Auslander algebra of $K[X]/(X^3)$. This has global dimension $2$ and is representation-finite. Then look at the Frobenius exact substructure described in \cite{ASo-Gor} on the Auslander algebra of a self-injective algebra.         
\end{exa}

\section{Parametrization using the Ziegler spectrum}

From now on, we impose the condition that $\mcC$ has weak cokernels, this means that for every morphism $f\colon X\to Y$ in $\mcC$ there exists a morphism $g \colon Y \to Z$ such that the following sequence is exact (in the middle) in the abelian category $\Modu \mcC^{op}$ (all covariant, additive functors $\mcC\to (Ab)$) 
\[
\Hom_{\mcC} (Z,-) \xrightarrow{\Hom (g,-)}\Hom_{\mcC} (Y,-) \xrightarrow{\Hom (f,-)} \Hom_{\mcC} (X,-) 
\]
This condition is equivalent to $\modu_1\mcC^{op}$ being abelian, in which case it also has enough projectives. In \cite{CB}, it has been shown that $\mcC$ has weak cokernels if and only $\mcA$ has products. This has been used in  \cite{CB} and \cite{Sch} to embed $\mcA=\overrightarrow{\mcC}$ in a locally coherent abelian category called the \emph{purity category}. We are not going to explain this further (as we have not used it).

\begin{lem} In a locally finitely presented additive category $\mcA$ with products, all definable categories are strongly definable. 
\end{lem}

To see this lemma: If $\mcA$ has products then every definable subcategory has products and 
${\rm fp}(\mcA)=\mcC$ has weak cokernel and by Lemma \ref{firstProperties}, every definable subcategory is a pure Serre subcategory of $\mcA$. 
To see that it is closed under taking pure-injective hulls we follow \cite[chapter 18, introduction]{Pr}: In this situation a definable subcategory has an intrinsic notion of purity, independent of how we realize it as a definable subcategory in a locally coherent abelian category. Furthermore, it is shown by \cite{Kr-Exactly}, in this situation, every definable category is \emph{exactly definable}, i.e.  equivalent to the category of exact functors on a small abelian category to abelian groups. In particular by Lemma \ref{PI-envelopes} this category is closed under pure-injective envelopes. 

\begin{rem}
Let us call a category \emph{really exactly definable} if it is equivalent to a category of exact functors from a small exact category to abelian groups. 
This leads to the open question: Is every strongly definable subcategory really exactly definable?  
\end{rem}

\begin{thm} (combine \cite[section (3.5), Lem 1]{CB} with \cite[Lem. 12.1.12]{Kr-book})
If $\mcC$ has weak cokernels then 
Ziegler-closed subsets form a the closed sets in a  topology on $\Zg(\mcA)$. This topological space is called the \textbf{Ziegler spectrum} of $\mcA$.
\end{thm}

We want to understand the exact substructures in cases where the Ziegler spectrum is known. For this we need the following: 
Let $\La$ be a ring, then we define the (left) Ziegler spectrum of $\La $ as $\Zg_{\La}:=\Zg (\La\Modu)$.  
\subsection{Examples}
\begin{exa}
As a consequence of \cite[Cor. 5.3.36, Cor. 5.3.37, Thm 5.1.12]{Pr} one obtains: For a finite-dimensional algebra $\La$ the following are equivalent
\begin{itemize}
    \item[(1)] $\La$ is of finite representation-type
    \item[(2)] $\Zg_{\La}$ is a finite set
    \item[(3)] $\Zg_{\La}$ does not contain any infinite-dimensional modules. 
\end{itemize}
In this case, $\Zg_{\La}$ is a discrete topological space, $\mcU_{max}$ consists of the indecomposable injectives in $\La\modu$. 
So, Ziegler-closed subsets containing $\mcU_{max}$ are in bijection with basic cogenerators in $\La\modu$. This is easily seen to be an equivalent description to Enomoto's theorem in this case. 
\end{exa}
\begin{exa} Ziegler spectrum in tame hereditary case has been described by Ringel in \cite{Ri-Ziegler}, we just look here at the easiest case: \\
We define $Q$ to be the Kronecker quiver $\xymatrix{1 &2\ar@<-3.5pt>[l]\ar[l]}$ and $\La=KQ$ for some infinite field $K$. Its Auslander-Reiten quiver (see picture below) has as vertices the indecomposables in $\La\modu$, they are divided in three types 1) $\mcP$ preprojectives (in the $\tau^-$-orbit of the projectives), they are denoted by their dimension vector $(n+1,n)$ , 2) $\mcR$ regulars, they are determined by the regular simple which they contain and their dimension vector, the regular simples are denoted by $S_{\lambda}, \lambda \in K\cup \{\infty\}=:\Omega$ , 3) $\mcQ$ preinjectives (in the $\tau$-orbit of the injectives), they are denoted by their dimension vector $(n,n+1)$. 

\[
\begin{tikzcd}[row sep=2em,column sep=0.001em]
&  &  &  &  &  & \vdots \arrow[d]  &\vdots\arrow[d]  & &  &  &  & & &   \\
& & & & &  & {S_{\lambda}[3]} \arrow[d] \arrow[u, shift right] \arrow[dashed, loop, distance=2em, in=305, out=235] & {S_{\mu}[3]} \arrow[d] \arrow[u, shift right] \arrow[dashed, loop, distance=2em, in=305, out=235] &  & &  &  &  & & \\
& {(2,1)} \arrow[rd, shift right] \arrow[rd] & & {(4,3)} \arrow[rd, shift right] \arrow[rd] \arrow[ll, dashed] & \cdots                    &  & {S_{\lambda}[2]} \arrow[d] \arrow[u, shift right] \arrow[dashed, loop, distance=2em, in=305, out=235] & {S_{\mu}[2]} \arrow[d] \arrow[u, shift right] \arrow[dashed, loop, distance=2em, in=305, out=235] &  &  & \cdots & {(3,4)} \arrow[rd, shift right] \arrow[rd]  &  & {(1,2)} \arrow[rd, shift right] \arrow[rd] \arrow[ll, dashed] & \\
{(1,0)} \arrow[ru, shift right] \arrow[ru] & & {(3,2)} \arrow[ru] \arrow[ru, shift left] \arrow[ll, dashed] &  & \cdots \arrow[ll, dashed] &  & S_{\lambda} \arrow[u, shift right] \arrow[dashed, loop, distance=2em, in=305, out=235] & S_{\mu} \arrow[u, shift right] \arrow[dashed, loop, distance=2em, in=305, out=235] & \cdots &  & \cdots \arrow[ru, shift right] \arrow[ru]   &  & {(2,3)} \arrow[ru, shift right] \arrow[ru]  \arrow[ll, dashed] &  & {(0,1)}\arrow[ll, dashed] \\
&  \mcP & & & &  & \mcR  & & & &  & \mcQ & & &        
\end{tikzcd}
\]
The arrows between the vertices indicate irreducible maps between the indecomposables and the dotted arrow the Auslander-Reiten translate, for every dotted arrow there is an almost split sequence. For more details look into \cite{Ri-tame}. \\
Then $\Zg_{\La}$ consists of the following points 
\begin{itemize}
\item[(1)] indecomposables in $\La\modu$
\item[(2)] For every $\lambda \in \Omega$ a Pr\"ufer module $S_{\lambda}[\infty]$, which is the filtered colimit (union) over $S_{\lambda} \infl S_{\lambda}[2]\infl S_{\lambda} [3] \infl \cdots$ 
\item[(3)] For every $\lambda \in \Omega$ an adic module $\hat{S_{\lambda}}$, which is the limit over $\cdots \defl S_{\lambda} [3] \defl S_{\lambda}[2]\defl S_{\lambda}$
\item[(4)] The generic module $G$, it is characterized by being an indecomposable module with $\Hom_{\La} (G, S_{\lambda})=0=\Hom_{\La} (S_{\lambda},G)$ for all $\lambda \in\Omega$
\end{itemize}
Now, $\mcU_{max}=\{ (1,2), (0,1)\}$ only consists of the two indecomposable injectives in $\La\modu$. Given a subset $\mcU\subseteq \Zg_{\La}$ containing $\mcU_{max}$ we find $T,M\subseteq \Omega $ 
\[
\begin{aligned}
\mcU^{fin} &:=\{U \in \mcU\mid U \in \La \modu, U \notin \mcU_{max}\}\\
\mcU_{T,M} &:=\mcU_{max}\cup \{S_{t}[\infty] \mid t\in T\}\cup \{ \hat{S_m} \mid m \in M\} \cup \{G\}
\end{aligned}
\]
such that $\mcU=\mcU_{max} \cup \mcU^{fin} \cup \mcU_{T,M}$ or 
$\mcU=\mcU_{max} \cup \mcU^{fin} \cup \mcU_{T,M}\setminus \{G\}$. 
Following Ringel's characterization in \cite{Ri-Ziegler} we find that $\mcU$ is Ziegler-closed iff   
\begin{itemize}
    \item[(a)] $\mcU^{fin}$ finite, then it and $T,M$ can be arbitarily  chosen (also empty is allowed) and 
    $\mcU=\mcU_{\max}\cup \mcU^{fin}\cup \mcU_{T,M}$ or if $T=M=\emptyset$ we can also have $\mcU= \mcU_{\max}\cup \mcU^{fin}$. 
    \item[(b)] $\mcU^{fin}$ infinite, then always $G \in \mcU$ but $T,M$ must satisfy the following. 
    \begin{itemize}
    \item[(c1)] If $\mcU^{fin}\cap \mcP$ is infinite, then $M=\Omega$ (all adics in) 
    \item[(c2)] If $\mcU^{fin}\cap \mcQ$ is infinite, then $T=\Omega$ (all Pr\"ufer in)
    \item[(c3)] For every $\lambda \in \Omega$, if $\mcU^{fin}\cap \{S_{\lambda} [n]\mid n \in \mathbb{N}\}$ is infinite, then $\lambda \in T\cap M$
    \end{itemize}
\end{itemize}
Before we start we need to understand some properties of the functors $\Hom_{\La} (-,U)$ and $\Ext^1_{\La}(-,U)$ for points of type (2,),(3),(4). In \cite[p. 46]{Ri-Inf} and \cite[section 3]{CB-Inf} we found the following vanishing where we set $S=S_{\lambda}$ and denote by $\mcR_{\lambda}$ be (the tube of) all regular modules with $S$ as a submodule. 
\[
\begin{aligned}
\Hom_{\La}(\mcR, G)=\Hom_{\La}(\mcQ, G) &=0 = \Ext^1_{\La}(\mcR ,G)=\Ext^1_{\La} (\mcP, G)\\
\Hom_{\La} (\mcR, \hat{S}) =\Hom_{\La}(\mcQ, \hat{S}) &=0=\Ext^1_{\La} (\mcP, \hat{S})= \Ext^1_{\La} (\mcR_{\mu}, \hat{S}) \quad \mu \neq \lambda \\ 
\Hom_{\La}(\mcR_{\mu}, S[\infty])=\Hom_{\La}(\mcQ, S[\infty]) &=0 = \Ext^1_{\La}(\mcR ,S[\infty])=\Ext^1_{\La} (\mcP, S[\infty]) \quad \mu \neq \lambda 
\end{aligned}
\]
As a consequence we see: 
If $U\in \{G, \hat{S}, S[\infty]\}$ and $\sigma$ a short exact sequence in $\La\modu$ with all three terms in either $\add(\mcP)$, $\add(\mcR)$ or $\add(\mcQ)$, then $\Hom (\sigma , U)$ is exact. \\
As there are very many Ziegler-closed sets in this case, we focus on two types: 
\begin{itemize}
\item[(I)] Either $\mcU^{fin}=\emptyset$, these give exact structures containing all almost split sequences.
\item[(II)] $\mcU=\overline{\mcU^{fin}}$, these are so-called Auslander-Solberg exact structures. Here, this is still an Auslander-Reiten category, the almost split sequences are precisely the ones of $\La\modu$ not starting in $\mcU^{fin}$. 
\end{itemize}

Now, we look at the exact structures in these cases: 
\begin{itemize}
\item[(I)] 
For $\mcU=\overline{\{U\}}$, we set $\Ext^1_U:=\Ext^1_{\mcE^{\mcU}}$. \\
We start with the unique maximal not abelian exact structure in this case $\mcU=\{G\}$, then $\Ext^1_G(\mcX, \mcY)= \Ext^1_{\La}(\mcX, \mcY)$ for all $(\mcX, \mcY) \in \{\mcP, \mcR, \mcQ\}^2\setminus \{(\mcQ, \mcP)\}$ and $\Ext^1_G(\mcQ, \mcP)=0$ (for this we leave the proof out). The interesting thing is that this is an exact substructure of global dimension $\geq 2$ (probably $=2$), since the following exact sequence $\sigma$ is not zero in $\Ext^2_G(\mcQ, \mcP)$: Let $R$ be a regular module, take a projective $\La$-module resolution and an injective $\La$-module resolution of $R$ and concatenate to an exact sequence $\sigma$
\[
P_1\infl P_0 \to I_0 \defl I_1
\]
Observe, this implies for all not abelian exact structures of type (I) that $\Ext^1( \mcQ, \mcP)=0$. \\
Next, we consider $\mcU=\overline{\{\hat{S}\}}$, we have $\Ext^1_{\hat{S}} (\mcX, \mcY) = \Ext^1_G (\mcX, \mcY)$ for all $(\mcX, \mcY)\neq (\mcR_{\lambda}, \mcP)$ and $\Ext^1_{\hat{S}} (\mcR_{\lambda}, \mcP) =0$.  \\
Now, we consider $\mcU=\overline{\{S[\infty]\}}$. Then $\Ext^1_{S[\infty]}(\mcQ, \mcR_{\lambda})=0$ and $\Ext^1_{S[\infty]}(\mcX, \mcY)=\Ext^1_{G}(\mcX, \mcY)$ for all $(\mcX, \mcY) \in \{ \mcP, \mcQ, \mcR_{\mu}, \mu \in \Omega\}^2\setminus \{(\mcP, \mcR_{\lambda})\}$. \\
In both cases $\mcU=\overline{\{\hat{S}\}}$ or $\mcU=\overline{\{S[\infty]\}}$ is the global dimension of $\mcE^{\mcU}$ still $\geq 2$. Just look at the exact sequence $\sigma$ as above. Choose in its definition $R$ to be a regular module $R$ with no summand in the tube $\lambda$, then for both exact substructures it gives an exact sequence which is not $2$-split. \\
Now, we look at intersections of these exact structures and respectively unions of the Ziegler-closed sets. \\
When $M=\Omega$ and $T=\emptyset$ then the exact structure consists of ses $\sigma = \sigma_p \oplus \sigma_{rq} $ such that $\sigma_p$ is an exact sequence in $\add (\mcP)$ and $\sigma_{rq}$ is an exact sequence in $\add (\mcR \cup \mcQ)$. It is very easily seen to be hereditary exact. \\
When $M=\emptyset$ and $T=\Omega$ then the exact structure consists of ses $\sigma = \sigma_{pr} \oplus \sigma_{q} $ such that $\sigma_{pr}$ is an exact sequence in $\add (\mcP\cup \mcR)$ and $\sigma_{q}$ is an exact sequence in $\add (\mcQ)$. It is very easily seen to be hereditary exact. \\
The case $T=M= \Omega$, then this is the minimal exact structure containing all almost split sequences. The short exact sequences in this structure are $\sigma_p \oplus \sigma_r \oplus \sigma_q$ with 
$\sigma_{p}$ is an exact sequence in $\add (\mcP)$, $\sigma_{r}$ is an exact sequence in $\add (\mcR)$ and $\sigma_{q}$ is an exact sequence in $\add (\mcQ)$. Again, we easily see that this is hereditary exact. 

\item[(II)] $\mcU^{fin}=:\mcH$. The exact structure corresponding to $\mcU$ is just given by all short exact sequences such that $\Hom(-,H)$ is exact on it for all $H$ in $\mcH$,
we write $\mcE=(\La \modu , F^\mcH)$. This case is well-studied in \cite{ASoI}, \cite{ASoII}, \cite{ASoIII}. If $\mcH$ is finite, then the exact structure always has enough projectives and enough injectives given by $\add (\mcH)$. Its global dimension can be characterized as follows $\gldim \mcE\leq k$ is equivalent to the following two conditions (i) $\gldim \End_{\La} (\bigoplus_{H\in \mcH} H) \leq k+2$ and (ii) $\idim_{\mcE} \La \leq k$. \\
For $\La=KQ$ with $Q$ the Kronecker quiver every global dimension can occur. This can be seen directly, just take 
$\mcH=\{(0,1),(1,2),(3,4), (5,6), (7,8), \cdots ,(2n-1, 2n)\}$. Then injective coresolutions are calculated via left $\add (\mcH)$-approximations and it can be easily seen that minimal injective coresolutions have at most $(n+1)$-injective modules, e.g. 
\[
(2n,2n+1) \infl (2n-1,2n)^{\oplus 2} \to (2n-3,2n-2)^{\oplus 2} \to\cdots \to (1,2)^{\oplus 2} \defl (0,1) 
\]
If you take $\mcH=\{S_{\lambda}, S_{\lambda}[3]\}$, then you find $\idim S_{\lambda} [2]= \infty$ and therefore we have infinite global dimension.\\
Another class of examples always gives hereditary exact substructures, take $\mcH=\{(n, n+1)\mid 0 \leq n \leq N\}$ for some $N\in \mathbb{N}$, then the cogenerator $\add(\mcH)$ is closed under quotients, this is easily seen to imply that the corresponding exact structure is hereditary exact. \\
Once you take $\mcH$ infinite, it is also easy to find infinite global dimensions: 
\[
\mcH= \{(2n-1,2n)\mid n \in \mathbb{N}\}
\]
Using minimal injective coresolutions for $(2n, 2n+1)$ for all $n \in \mathbb{N}$, we find objects of injective dimension $n$ for every $n \in\mathbb{N}$, this implies $\gldim =\infty$. 
\end{itemize}
\end{exa}

\begin{exa} 
We describe all exact substructures on finitely generated modules over a commutative discrete valuation ring $R$ with maximal ideal $P$.  We recall the description of the Ziegler spectrum from \cite[Section 5.2]{Pr}:\\
The points in $\Zg_R$ are:
\begin{itemize}
\item[(a)] indecomposable modules of finite length $R/P^n$, $n\geq 1$
\item[(b)] the $P$-adic completion $\overline{R}=\lim R/P^n$ (this is the limit over $\cdots \to R/P^{2} \to R/P=k$)
\item[(c)] The Pr\"ufer module $R_{P^{\infty}}=\colim R/P^n$ (this is the colimit over $k=R/P \to R/P^2 \to \cdots $)
\item[(d)] the quotient division ring $Q=Q(R)$ of $R$
\end{itemize}
Now, $\mcU_{max}=\{Q, R_{P^{\infty}}\}$ is the Ziegler-closed set given by the indecomposable injective $R$-modules. We also observe that $\Zg_R':=\{\overline{R}\}\cup \mcU_{max}$ is Ziegler-closed.\\
Next, all Ziegler-closed subsets containing $\mcU_{max}$ are given by: 
\begin{itemize}
\item[(1)] for $\emptyset \neq L\subseteq \mathbb{N}$ \emph{finite}, we have 
\[
\mcU_L:=\{R/P^n \mid n \in L\} \cup \mcU_{max}
\]
\item[(2)] for $\emptyset \subseteq L\subseteq \mathbb{N}$ \emph{arbitrary subset} we have 
\[
\mcV_L:=\{R/P^n \mid n \in L\} \cup \Zg_R' 
\]
\end{itemize}
So let us describe the exact structures on $\mcC=R\modu$ the category of finitely generated left $R$-modules corresponding to these closed sets: 
\begin{itemize}
    \item[(max)] Trivially $\mcU_{max}$ corresponds to the abelian structure on $\mcC$, this is hereditary and with enough projectives (but not with enough injectives) 
    \item[(min)] and $\Zg_R$ corresponds to the split exact structure. 
    \item[($\Zg$')] The Ziegler-closed set $\Zg_R'$ corresponds to the exact structure $\mcE'$ making the torsion functor exact. This is a hereditary exact structure , cp. example \ref{tpExa}.
    \item[($\mcU_L$)] This corresponds to the exact substructure $\mcE_L$ such that $\Hom_R(-,R/P^n)$, $n \in L$ are exact functors to abelian groups. 
    \item[($\mcV_L$)] This corresponds to the exact substructure $\mcE_L'$ such that the torsion functor and \\
    $\Hom_R(-,R/P^n)$, $n \in L$ are exact.  
\end{itemize}
First, of all in general how can one see that for $\emptyset\neq L \subseteq \mathbb{N}$ finite: $\mcE_L$ and $\mcE_L'$ are different exact structures?\\
Take a short exact sequence $R\infl R \defl R/P^n$ and the pushout along $R\defl R/P^m$, this gives an exact sequence $R/P^m \infl R/P^{m+n} \defl R/P^n$. But the cartesian square induces another exact sequence 
\[
R\infl R/P^m \oplus R \defl R/P^n
\]
It is easily seen to be not exact in $\mcE'$ if $n\neq m$ and we conclude that $\Ext^1_{\mcE'}(R/P^n,R)=0$. But if you apply 
$\Hom_R(-,R/P^{\ell})$ using $\Hom_R(R/P^s, R/P^{\ell})=R/P^{\min (s,\ell )}$ we see that this is exact for $n,m$ both larger or equal than $\ell$. \\
We say that $L$ has \emph{gaps} if there is an interval $[a,b]$ such that $[a,b]\cap L=\emptyset$ and $b+1\in L$ and $a-1\in L$ if $a>1$. 
If $L$ has gaps then $\gldim \mcE_L=\infty=\gldim\mcE_L'$ (for $\mcE_L'$ we also allow $L$ to be an infinite subset).\\
In this case one can always find an infinite injective $\mcE_L^{(')}$-coresolution for an $R/P^s$ some $s\in [a,b]$. We give them as sequence of short exact sequences.  
\begin{itemize}
\item[($a=1$)] Take $s=b$ and $R/P^b\infl R/P^{b+1} \defl R/P$, then continue with 
$R/P \infl R/P^{b+1} \defl R/P^b$ and repeat with the first short exact sequence etc. 
\item[(2)] If $a>1$ and $b+a$ even then we take $s=\frac{1}{2} (a+b)$ and the short exact sequence 
\[R/P^s\to R/P^{a-1}\oplus R/P^{b+1}\defl R/P^s\]
and then continue with the same sequence. 
\item[(3)] If $a>1$ and $b+a$ uneven then we take $s=\frac{1}{2}(b+a-1)$ and first $R/P^s\infl R/P^{a-1}\oplus R/P^{b+1} \defl R/P^{s+1} $ then 
$R/P^{s+1} \infl R/P^{a-1} \oplus R/P^{b+1} \defl R/P^s$ and then repeat with the first sequence. 
\end{itemize}

If $L$ has no gaps and $L \neq \mathbb{N}$ then $L=[1,n]$ for some $n\in \mathbb{N}$. The set $\{R/P^{\ell}\mid \ell\in L\}$ is closed under quotients, so injective coresolutions in this class of modules will always end after one short exact sequence. We show in the next Lemma that these exact structures are hereditary exact. 
\end{exa}

\begin{lem} \label{radExt} Let $R$ denote a commutative discrete valuation ring   with maximal ideal $P$ and let $n \in\mathbb{N}$.  
We have a functor $\rad^n \colon R\modu \to R\modu$ defined by $\rad^n M = P^n M$. Let $L=[1,n]\subseteq \mathbb{N}$. 
\begin{itemize}
\item[(1)] Then we have $\Ext^1_{\mcE_L} (-,-)=\rad^n\Ext^1_R(-,-)$ and $\Ext^1_{\mcE_L'} (-,-)=\rad^n_R\Ext^1_{\mcE'}(-,-)$
\item[(2)] For every $\mcE_L$-exact sequence $\sigma$ and every object $X$, the sequences $\Ext^1_{\mcE_L} (X,\sigma)$ is right exact.\\
\noindent
For every $\mcE'_L$-exact sequence $\sigma$ and every object $X$, the sequences $\Ext^1_{\mcE'_L} (X,\sigma)$ is right exact.
\end{itemize}    
In particular, $\mcE_L$ and $\mcE'_L$ are hereditary exact. 
\end{lem}

\begin{proof}
\begin{itemize}
\item[(1)] We first describe the $R$-module structure on 
\begin{itemize}
    \item[(a)] $\Ext^1_R (R/P^m, R/P^{\ell})\cong R/P^s$ where $s=\min (m, \ell)$. For $a,b\in \mathbb{N}$ we write $\sigma_{a,b}\in \Ext^1_R(R/P^m, R/P^{\ell})$ for an exact sequence with middle term $R/P^a\oplus R/P^b$ (whenever this exist). In $R/P^s$ we pick $P=(p)$ and we have the following mult. by $p$ \\
    $1 \mapsto p \mapsto p^2 \mapsto \cdots \mapsto p^{s-1} \mapsto 0$
    this corresponds to the following on the Ext-group
    \begin{itemize}
     \item[(a1)] If $s=\ell \leq m$ we have \\
      $\sigma_{0, m+\ell}\mapsto \sigma_{1, m+\ell-1} \mapsto\cdots \mapsto  \sigma_{\ell-1, m+1}  \mapsto 0$\\
     Then we have ${\rm rad}^n\Ext^1_R(R/P^m, R/P^{\ell})\cong R/P^{s-n}$ whenever $n<s$ and zero otherwise. For $n<s=\ell$ it is the image of $p^n$, i.e.   \\
      $\sigma_{n, m+\ell-n}\mapsto \sigma_{n+1, m+\ell-n-1} \mapsto\cdots \mapsto  \sigma_{\ell-1, m+1} $
     \item[(a2)] If $\ell >m=s$ we have \\
     $\sigma_{\ell +m, 0}\mapsto \sigma_{\ell+m-1, 1} \mapsto\cdots \mapsto  \sigma_{\ell+1,m-1} \mapsto 0$\\
     Then we have ${\rm rad}^n\Ext^1_R(R/P^m, R/P^{\ell})\cong R/P^{s-n}$ whenever $n<s$ and zero otherwise. For $n<s=\ell$ it contains the following elements  \\
      $\sigma_{\ell+m-n, n}\mapsto \sigma_{\ell +m-n-1,n+1} \mapsto\cdots \mapsto  \sigma_{\ell+1, m-1} $
    \end{itemize}
    \item[(b) ] $\Ext^1_R(R/P^m, R)\cong R/P^m$ we write $\sigma_a$ for the extension with $R \oplus R/P^a$ as middle term. Then the multiplication by $p$ is given by \\
    $\sigma_{0} \mapsto \sigma_{1} \mapsto \cdots \mapsto \sigma_{m-1} \mapsto 0$\\
    The submodule ${\rm rad}^n\Ext^1_R(R/P^m, R)$ is of course zero is $n\geq m$ and if $n<m$ is given by the following \\
    $\sigma_{n} \mapsto \sigma_{n+1} \mapsto \cdots \mapsto \sigma_{m-1}$
\end{itemize} 
We claim $\Ext^1_{\mcE_{[0,n]}}=\rad^n\Ext^1_R$. First observe that in $\mcE_L$: $R/P^a$ $1\leq a \leq n$ are injectives and they are also projectives. So for $s= \min (\ell ,m) \leq n$ we have $\Ext^1_{\mcE_{[0,n]}}(R/P^m, R/P^{\ell})=0$ and for $m\leq n $ we have $\Ext^1_{\mcE_{[0,n]}}(R/P^m, R)=0$. \\
So we may always assume wlog that $n< \min (\ell,m)$, then proceed by induction over $n$. For $n=1$, a short exact sequence of in $R\modu$ is in $\mcE_{[0,1]}$ iff the indcomposable summands of number of the indec. summands of the outer terms add up to the indec summands of the middle term. This means in case (a) all exact sequence are in this exact structure except $\sigma_{0, m+n}$ and $\sigma_{m+n,0}$, in case (b) all except $\sigma_0$. \\
For $n>1$, it is enough to observe that $\Hom (-, R/P^n)$ is 
\begin{itemize}
\item[(ad a1)] exact on $\sigma_{n, m+\ell -n}$ and not exact on  $\sigma_{n-1, m+\ell -n+1}$
\item[(ad a2)] exact on $\sigma_{\ell+m-n, n}$ and not exact on  $\sigma_{\ell +m -n+1, n-1}$
\item[(ad b)] exact on $\sigma_{n}$ and not exact on  $\sigma_{n-1}$
\end{itemize}
Then the rest follows by induction hypothesis. \\
Now, for $\mcE'_L$ one observes that $\Ext^1_{\mcE'_L}(X,Y)=\Ext^1_{\mcE_L}(X,Y)$ for all $X,Y$ torsion, and for $Y$ free it is zero.\\
As $\Ext^1_{\mcE'}(X,Y)=\Ext^1_R(X,Y)$ for all $X,Y$ torsion and for $Y$ free it is zero. Then the claim $\Ext^1_{\mcE'_L}$ follows from the proof for $\Ext^1_{\mcE_L}$. 
\item[(2)] Taking $n$-th radical of an epimorphism in $R \modu$ is again an epimorphism - as the $n$-th radical can be described as the image of multiplication by $p^n$ (making it also a quotient and not ony a submodule). As $\Ext^1_R(X, \sigma)$ (resp. $\Ext^1_R(\sigma, X)$) is right exact for all exact sequences $\sigma$, we conclude that $\rad^n \Ext^1_R(X, f)$ (resp. $\rad^n \Ext^1_R(g, X)$) are epimorphisms for $f$ an epimorphism and $g$ a monomorphism. \\
Of  course, on general exact sequence $\rad^n$ is not a middle-exact functor (but this follows from (1) since for $\sigma \in \mcE_L$ the sequences $\Ext^1_{\mcE_L}(X,\sigma)$ and $\Ext^1_{\mcE_L}(\sigma ,X)$ are middle exact). In particular, the right exactness claim follows. \\
For $\Ext^1_{\mcE'_L}(X,\sigma) $ we can restrict to $X$ indecomposable torsion and as $\sigma$ in $\mcE'$, it fits into an exact sequence of ses $\sigma_{tor}\infl \sigma \defl \sigma_{free}$ with $\sigma_{tor}$ the torsion part and $\sigma_{free}$ the free part, and we conclude that $\Ext^1_{\mcE'_L}(X,\sigma)= \Ext^1_{\mcE'_L} (X, \sigma_{tor}) =\Ext^1_{\mcE_L} (X, \sigma_{tor})$ is right exact.   
\end{itemize}
\end{proof}

After understanding exact substructures of $R\modu$ for $R$ a commutative discrete valuation ring, we are ready to generalize this to commutative Dedekind domains:  
\begin{exa}
We describe all exact structures on finitely generated left modules over a commutative Dedekind domain. \\
The Ziegler spectrum is studied as a more general case of the discrete valuation ring, again we follow \cite[Section 5.2]{Pr} for its description: Let $\mSpec (R):=\{P\mid \text{ max. ideal in } R\}$. 
The points in $\Zg_R$ are:
\begin{itemize}
\item[(a)] indecomposable modules of finite length $R/P^n$, $n\geq 1$, and $P\in \mSpec (R)$.
\item[(b)] the $P$-adic completion $\overline{R}_P=\lim R/P^n$ (this is the limit over $\cdots \to R/P^{2} \to R/P$) for $P\in \mSpec (R)$.
\item[(c)] The Pr\"ufer module $R_{P^{\infty}}=\colim R/P^n$ (this is the colimit over $R/P \to R/P^2 \to \cdots $) for $P\in \mSpec (R)$.
\item[(d)] the quotient division ring $Q=Q(R)$ of $R$
\end{itemize}
Now, $\mcU_{max}=\{Q\}\cup  \{R_{P^{\infty}}\mid  P\in \mSpec (R)\}$ is the Ziegler-closed set given by the indecomposable injective $R$-modules. We describe all Ziegler-closed subsets containing $\mcU_{max}$ (following loc. cit.). \\
First we fix some notation, let $L\subseteq {\rm mSpec} (R) \times \mathbb{N}$ always denote such a subset and for $P\in {\rm mSpec} (R)$ let $L_P:=\{\ell \in \mathbb{N}\mid (\ell, P)\in L\}$. 
Subsets of indecomposable finite length modules are of the form $\mcF_L=\{ R/P^{\ell} \mid (P, \ell) \in L\}$.  
We fix a closed subset $\mcU$ and denote by $\mcF_L$ its points of finite length. 
\\
(type 1) $L=\emptyset$. 
For every $M\subseteq \mSpec  (R)$ we have a closed subset $\Zg_M'=\mcU_{max} \cup \{\overline{R}_P \mid P \in M\}$. We define  $\Zg':=\Zg'_{\mSpec (R)}$. \\
(type 2) $0< \lvert L \rvert <\infty$. 
Then $\mcU=\mcF_L \cup \Zg'_M$ for an arbitrary subset $M \subseteq \mSpec (R)$ (the sets $L$ and $M$ are independent from each other).\\
(type 3) $\lvert L \rvert =\infty$. 
We define $M_L:=\{P\in \mSpec (R)\mid \exists \; n \in \mathbb{N} \colon (P,n)\in L \}$ and in this case 
$\mcU=\mcF_L\cup \Zg'_{M_L}$.\\

Let us look at the corresponding exact substructures: 
\begin{itemize}
\item[(type 1)] 
Make all $P$-torsion functors exact for all $P\in M$. As we are dealing with hereditary torsion pairs, the torsion functors are left exact and example \ref{tpExa} applies to show that these are hereditary exact structures. 
\item[(gaps)]
Let us assume we are in type 2 or type 3.\\
If for some $P\in \mSpec (R)$ we have $L_P$ has a gap (see previous example), then we find an infinite injective coresolution as in the previous example and it follows $\gldim = \infty$.   
\item[(no gaps)]
Let us assume we are in type 2 or type 3. If all non-empty $L_P$ have no gaps for every maximal ideal $P$, then we find $L_P=[1, n_P]$ for an $n_P\in \mathbb{N}$. We claim that we only have hereditary exact structures in this case. We give the proof in the next Lemma.

\end{itemize}
\end{exa}

Now, let $R$ be a commutative Dedekind domain. Observe that $\Ext^1_R$ only takes values in torsion modules. We write $()_P$ for its $P$-torsion submodule and $()_{tor}$ for the remaining torsion summands. So, for two finitely generated $R$-modules $X=R^t\oplus X_P \oplus X_{tor},Y=R^s\oplus Y_P \oplus Y_{tor}$ we have 
\[
\begin{aligned}
\Ext^1_R(X,Y) &=\Ext^1_R(X,Y)_P\oplus \Ext^1_R(X,Y)_{tor}\\
\Ext^1_R(X, Y)_P &= \Ext^1_R(X_P, Y_P\oplus R^s)\\
\Ext^1_R(X,Y)_{tor} &=\Ext^1_R (X_{tor}, Y_{tor}\oplus R^s)
\end{aligned}
\]
Then $\Ext^1_R = (\Ext^1_R)_P\oplus (\Ext^1_R)_{tor}$ is a direct sum decomposition of bifunctors (but this does not imply that these subfunctors are middle exact for the abelian structure on $R\modu$). Let $F_P\colon R\modu\to R\modu$ be the functor $X\mapsto R^t \oplus X_P$ (resp. $F_{tor}(X)=R^t\oplus X_P$), induced by the projection into the torsionfree part (in the split hereditary torsion pairs considered in (type 1)). They preserve epimorphisms. In particular, if $f\colon A\defl B$ is an epimorphism, then the following are also epimorphism as $\Ext^1_R(M,-)$ preserves epimorphisms for all objects $M$ 
\[
\begin{aligned}
\Ext^1_R(X, f)_P &=\Ext^1_R(F_P(X), F_P(f))\\
\Ext^1_R(X,f)_{tor} &= \Ext^1_R(F_{tor}(X), F_{tor} (f))
\end{aligned}
\]

As it looks simpler let us look at the case of only one prime:
\begin{lem} \label{P-rad} Let $R$ be a commutative Dedekind domain. Let $P$ be a fixed maximal prime ideal. We define for $M \in R\modu, P \in \mSpec (R), n \in \mathbb{N}$ the following $\rad_P^n M:= P^nM$. If $L=\{P\} \times [1,n]$ we denote by $\mcE_L$ the exact structure corresponding to $\mcF_L \cup \mcU_{max}$. 
\begin{itemize}
\item[(1)]    
Then $\Ext^1_{\mcE_L}=(\rad_P^n \Ext^1_R)\oplus (\Ext^1_R)_{tor}$
\item[(2)] The exact structure $\mcE_L$ is hereditary exact. 
\end{itemize}
\end{lem}

\begin{proof}
\begin{itemize}
\item[(1)]
As for different primes the torsion submodules are $\Hom_R$- and $\Ext_R$-orthogonal, the exactness of $\Hom(-,R/P^a)$ for some $a\in \mathbb{N}$ only depends on the $P$-torsion and the free module summand. The same proof as in Lemma \ref{radExt} applies.    
\item[(2)] By the discussion before the Lemma and knowing that $\rad^n$ preserves epimorphisms, it follows from (1) that $\Ext^1_{\mcE_L}$ preserves epimorphisms. Therefore $\mcE_L$ is hereditary exact. 
\end{itemize}
\end{proof}

But it actually is the same for an arbitrary subset of primes: 
\begin{lem} \label{herSubstr} Let $R$ be a commutative Dedekind domain. \\
Let $L=\bigcup_{P\in M} \{P\} \times [1,n_P]\subseteq \mSpec (R) \times \mathbb{N}$  and $\mcE=\mcE^{\mcU}$ for a Ziegler-closed subset $\mcU$ with modules of finite length given by $\mcF_L$, then $\mcE$ is hereditary exact. 
\end{lem}
\begin{proof} (of Cor. \ref{herSubstr}) $\mcU=\mcF_L \cup \Zg'_M$ for some subset $M\subseteq \mSpec (R)$ (type 2 or type 3). \\
Let us denote by an index $tor(M^c)$ the torsion summand corresponding to the complement of $M$. The above Lemma generalizes to 
\[\Ext^1_{\mcE_L}=\bigoplus_{P \in M} (\rad_P^{n_P} \Ext^1_R)\oplus (\Ext^1_R)_{tor (M^c)}\]
This is clear as intersection of exact substructures correspond to intersecting the corresponding $\Ext^1$-subfunctors. \\
This functor is still preserving epimorphisms as before. 
\end{proof}

\subsection{Exact structures making additive functors exact}


We would like to introduce a construction of an exact structure following \cite{DRSS}.\\
\noindent
Let $(\mcA, \mcS)$ and $(\mcB, \mcT)$ be idempotent complete exact categories and $f\colon \mcA\to \mcB$ an additive functor. We denote by $\mcS_f=\mcS_{f, \mcT}\subseteq \mcS$ the class of exact sequences $\eta$ such that $f(\eta)$ is in $\mcT$. Observe, that this depends also on the exact structures and not just on the additive functor. 

\begin{lem} (\cite[Lem.1.9, Prop.1.10]{DRSS})
Then the following are equivalent
\begin{itemize}
    \item[(1)] $\mcS_f$ is an exact structure
    \item[(2)] Given a short exact sequence in $\mcS_f$ and a morphism to the third object of the sequence (resp. starting at the first object in the sequence), then the pullback (resp. pushout) of the short exact sequence is in $\mcS_f$.  
\end{itemize}
\end{lem}
It is also straight-forward to prove that (2) implies (1) by checking that compositions of $\mcS_f$-deflations (resp. -inflations) are $\mcS_f$-deflations (resp. -inflations).\\
We will now look at particular situations ensuring that $\mcS_f$ becomes an exact structure. 

\begin{cor} \label{split}
Assume that $f\colon (\mcA, \mcS)\to (\mcB, \mcT_0)$ is an exact functor between exact categories and let  $\mcT\subseteq \mcT_0$ be an exact substructure then $\mcS_f$ considered with respect to $\mcT$ is an 
 exact structure. 
\end{cor}
Given an exact category $\mcE$, we write $\rm{ex} (\mcE)$ the poset of exact substructures of $\mcE$. 
The previous corollary translates to: Every exact functor $f\colon \mcE\to \mcF$ induces a morphism of posets 
\[
f^* \colon \rm{ex} (\mcF) \to \rm{ex} (\mcE) , \quad \mcT \mapsto \mcS_{f, \mcT}
\]

\begin{proof}
We take $\eta \colon 0 \to A \to B \to C \to 0$ in $\mcS$ such that $f(\eta) \in \mcT$ and a morphism $c\colon C'\to C$. The pullback $\eta'$ of $\eta $ along $c$ exists and maps under $f$ to the pullback of $f(\eta) \in \mcT_0$ along $f(c)$. But $\mcT$ is an exact substructure of $\mcT_0$ and as $f(\eta) \in \mcT$ it follows that $f(\eta')\in \mcT$. The rest of the proof is the dual statement.  
\end{proof}

\begin{lem} \label{leftEx}
Let $f\colon (\mcA, \mcS) \to (\mcB, \mcT)$ be an additive functor between exact categories and assume that $\mcB$ is weakly idempotent complete. 
Assume either  
\begin{itemize}
    \item[(1)] $f$ is right exact (i.e. if $0\to A\to B \to C \to 0$ is in $\mcS$, then $f(A) \to f(B)$ has an image in $\mcB$ and $\Bild (f(A)\to f(B)) \infl f(B) \defl f(C) $ is in $\mcT$), or
    \item[(2)] $f$ is left exact (i.e. the opposite functor between the opposite exact categories is right exact) 
\end{itemize}
Then $\mcS_f$ is an exact structure.
\end{lem}

\begin{proof}
We prove only (1) as (2) is analogous. 
Assume $f$ is right exact. 
Take $\eta \colon A \infl B \defl C $ in $\mcS_f$. Let $\gamma \colon C'\to C$ be any morphism. Now, we pull-back $\eta \in \mcS$ along $\gamma $ to an $\eta' \colon A \infl B' \defl C' $ in $\mcS$. Applying $f$ gives a commutative diagram with exact rows 
\[ 
\begin{tikzcd}
f(A) \arrow[r, "a"] \arrow[d,-,double equal sign distance,double]  & f(B') \arrow[r, two heads] \arrow[d, "b"'] & f(C') \arrow[d] \\
f(A) \arrow[r, tail]          & f(B) \arrow[r, two heads]                  & f(C)           
\end{tikzcd}
\]
Since $b\circ a$ is an inflation it follows that $a$ is an inflation (using that $\mcB$ is idempotent complete) and therefore $\eta'\in \mcS_f$.
Now, we pushout $\eta $ along a morphism $\alpha \colon A\to A''$ to an $\eta''\colon A'' \infl B'' \defl C $ in $\mcS$. Recall from \cite[Prop. 2.12]{Bue},  that we also get an an induced short exact sequence $\tilde{\eta}\colon  A \infl A''\oplus B \defl B'' $. \\
Applying $f$ gives a commutative diagram with exact rows
\[ 
\begin{tikzcd}
f(A) \arrow[r, "d", tail] \arrow[d, "c"] & f(B) \arrow[r, two heads] \arrow[d] & f(C) \arrow[d,-,double equal sign distance,double] \\
f(A'') \arrow[r]                         & f(B'' \arrow[r, two heads]          & f(C)          
\end{tikzcd}
\]
Let $D$ be the push-out of $d,c$, since $\mcT$ is an exact structure, we get an exact sequence $ f(A) \infl f(A'') \oplus f(B) \defl D $ (in $\mcT$). In particular, the first map is an inflation. Now apply $f$ to $\tilde{\eta}$ to get a right exact sequence $f(A) \to f(A'') \oplus f(B'') \defl f(B'') $, since the first map is an inflation, we conclude that $\tilde{\eta} \in \mcS_f$ and $D\cong f(B'')$. This implies that that the lower row coincides with the push-out short exact sequence in $\mcT$, in particular $\eta'' \in \mcS_f$.

\end{proof}

\begin{exa} \label{tpExa} 
Given an abelian category $\mcB$, Dickson in \cite{Dic} defined a torsion pair to be a pair of two full subcategories $(\mcT, \mcF)$ satisfying $\Hom (\mcT, \mcF)=0$ and for every $B$ in $\mcB$ there exists an exact sequence $B_t\infl B \defl B_f$ with $B_t$ in $\mcT$, $B_f$ in $\mcF$. Then $\mcT$ and $\mcF$ are fully exact subcategories in $\mcB$ but they may or may not be homologically exact. 
The assignment $B \mapsto B_t$ extends to a functor $t \colon \mcB \to \mcT$ (cf. loc. cit. Cor. 2.5), this functor preserves inflations but is in general not left exact. But it is left exact for a so-called \emph{hereditary torsion pair} (cf.\cite[ Ch. VI, Prop. 1.7]{St}) - characterized by $\mcT$ being closed under subobjects. So assume we have an hereditary torsion pair.    
By Lemma \ref{leftEx}, we have an exact substructure $\mcB_t=(\mcB, \mcS_t)$. In $\mcB$ every exact sequence $X\infl Y \defl Z$ gives rise to a left exact sequence $t(X) \infl t(Y) \to t(Z)$ and a right exact sequence $X/t(X) \to Y/t(Y) \defl Z/t(Z)$. Short exact sequences in $\mcB_t$ are those for which both these sequences are short exact sequences (using the snake Lemma), i.e. they can be seen as an extension between an exact sequence in $\mcT$ and one in $\mcF$. By definition, both inclusion functor $\mcT \to \mcB_t$ and $\mcF \to \mcB_t$ are exact and the pair of subcategories $(\mcT, \mcF)$ in $\mcB_t$ is a torsion pair in an exact category (cf. e.g. \cite[Ex. 3.5]{AET}). Now, we show: 
\[
\gldim \mcB_t =  {\rm max} (\gldim \mcT, \gldim \mcF )
\]
As $t\colon \mcB_t \to \mcT$ is an exact functor with $t\circ i\cong \id_{\mcT}$ it follows that $\mcT$ is homologically faithful in $\mcB_t$ (and similarly $\mcF$ is homologically faithful). But they are also homologically exact in $\mcB_t$ as the map on $\Ext^n$ is also obviously surjective (given $[\eta] \in \Ext^n_{\mcB_t} (t(X), t(Y))$, then $[t(\eta)]$ is in the image $\Ext^n_{\mcT} (t(X), t(Y))\to \Ext^n_{\mcB_t} (t(X), t(Y))$ but by definition $t(\eta)$ and $\eta$ are equivalent because we have a morphism $t(\eta) \to \eta$ of exact sequences with fixed end terms). 
This implies $\gldim \mcB_t\geq {\rm max} (\gldim \mcT, \gldim \mcF )$. We claim that we have equality. To this, we define 
$n-1={\rm max} (\gldim \mcT, \gldim \mcF )<\infty $. So, given an $n$-exact sequence in $\mcB_t$, say $\eta \colon X \infl X_1 \to \cdots \to X_n \defl Y$, we want to see that $[\eta]=0$ in $\Ext^n_{\mcB_t}(Y,X)$. We have $t(\eta)$ and $\eta/t(\eta)$ are both exact sequences which are composed of short exact sequences in $\mcT$ (resp. $\mcF$), so they are in the image of $\Ext^n_{\mcT} (t(Y),t(X))\to \Ext^n_{\mcB_t} (t(Y),t(X))$, but as $\Ext^n_{\mcT}=0=\Ext^n_{\mcF}$ it follows that both $[t(\eta)]=0$ and $[\eta/t(\eta)]=0$, let us call this property $n$-split. But $n$-split exact sequences are closed under extensions (this is an easy exercise using the derived category $D^b(\mcB_t)$). It follows that $[\eta]=0$. \\
We apply this in the following situation: Given $R$ a principal ideal domain  $\mcB=R\modu$ and $\mcT$ be the full subcategory of torsion modules (i.e. they have a non-zero annihilator) and $\mcF$ the subcategory of free modules. This gives a hereditary torsion pair in $\mcB$. We conclude that in this case $\gldim \mcB_t=1$, i.e. $\mcB_t$ is still an hereditary exact category. 
\end{exa}

\section{Appendix: Ind-Completion of (small) exact categories}

These notes are based on the recent preprint of Positselski \cite{P-locCoh} - but we prefer the construction using the Gabriel-Quillen embedding (this way, we extend Crawley-Boeveys classical dictionary to exact structures \cite{CB}). 

\subsection*{Locally finitely presented additive categories. }

Here, we give a quick summary of the \emph{correspondence} from \cite{CB}.\\

Let $\mcC$ an essentially small additive category. We define $\hat{\mcC}:=\Modu \mcC$ to be the category of all additive functors $\mcC^{op} \to (Ab)$, we call this the category of (left) $\mcC$-\textbf{modules}. It is easily seen to be an abelian category.We have the (covariant) Yoneda embedding 
\[
\mathbb{Y} \colon \mcC \to \hat{\mcC}, \quad X \mapsto (-,X):= \Hom_{\mcC} (-, X) 
\]
this is fully faithful, the essential image consists of (some) projective objects which we call \textbf{representable functors}.\\ 
 Every object in $\hat{\mcC}$ is as a small colimit of representables 
 - for $F \in \hat{\mcC}$ define the \emph{slice} category $\mcC / F$ for $F\in \hat{\mcC}$ (objects: $(X,x)$, $X\in \mcC , x\in F(X)$, morphisms $f\colon X\to X'$ in $\mcC$ such that $F(f)(x')=x$), then we have a small category and a functor $\Phi \colon \mcC/F\to \hat{\mcC}, (X,x)\mapsto (-,X)$. Its colimit is $F=\colim_{\mcC/F} \Phi$.

\begin{dfn} (\cite{SGA4}, Expose I)
We define the \textbf{ind-completion} $\overrightarrow{\mcC}$ (in the literature denoted as ${\rm Ind}(\mcC)$) as the closure of $\mcC$ under arbitrary directed colimits: 
Objects are functors $D\colon I \to \mcC$ from small filtered categories $I$. Morphisms are defined as 
\[
\begin{aligned}
\Hom(D\colon I \to \mcC, E\colon J \to \mcC) &:= \Hom_{\Modu \mcC} (\colim_I \mathbb{Y}D, \colim_J \mathbb{Y}E)\\ &= \lim_{i\in I} \colim_{j\in J} \Hom_{\mcC} (D(i) , E(j)) 
\end{aligned}
\]
Observe that the Yoneda embedding factors over $\overrightarrow{\mcC}$. 
Via the Yoneda embedding, we can identify this with the following full subcategory of $\hat{\mcC}$
\[
\overrightarrow{\mcC}:= \{\colim_{i\in I}(-,X_i)\mid \quad 
(X_i)_{i \in I}\; I\!-\!\text{shaped diagram in }\mcC \text{ with } I\text{ directed set }\}
\]    
\end{dfn}
\begin{rem}
The second description uses that closure under small filtered colimits it the same as closure under small directed colimits, cf. \cite{ARo}, Thm 1.5.    \end{rem}

\begin{pro} (\cite{SGA4}, Expose I, Prop. 8.6.4)\label{SGA4Ind}
Ind-completion is a 2-functorial. \\
An additive functor $f$ is faithful (resp. fully faithful) if and only if $\overrightarrow{f}$ is faithful (resp. fully faithful). \\
Furthermore the ind-completion $\overrightarrow{f}$ of an additive functor $f\colon \mcC\to \mcD$ is an equivalence if and only it is fully faithful and the essential image inclusion $\Bild f \to \mcD$ induces an equivalence on idempotent completions.  
\end{pro}

 Let $\mcD$ be an additive category, we denote by ${\rm Add}(\mcC, \mcD)$ the category of additive functors from $\mcC$ to $\mcD$ and ${\rm Add}_{fc}(\mcC, \mcD)$ the subcategory of functors which preserve directed colimits (whenever these exist). 

Ind-completion can be defined for arbitrary additive and even arbitary categories categories and can be characterized by a universal property such as:   
\begin{lem} (\textbf{Universal property of ind-completion}, cf. \cite{SGA4}, Expose I, Prop. 8.7.3) \label{UP-add}
Let $\mcC$ be a small additive category, then $\overrightarrow{\mcC}$ has all directed colimits. \\
Assume that $\mcD$ is an additive category which is closed under arbitrary directed colimits. Precomposition with $\mcC\to \overrightarrow{\mcC}$ is an equivalence of categories 
\[
{\rm Add}_{fc}(\overrightarrow{\mcC}, \mcD) \to {\rm Add}(\mcC, \mcD)
\]
\end{lem}
Furthermore, it also has the following property 
\begin{lem} (\cite{CB}, Lem. 1)\label{ic}
$\overrightarrow{\mcC}$ is idempotent complete. 
\end{lem}

For small additive categories, we have an alternative description of the ind-completion found in \cite{CB}.
\begin{dfn}
Let $\mcC$ be a small additive category. We say that an object $F$ in $\hat{\mcC}$ is \textbf{flat} if the tensor functor $F\otimes_{\mcC} -\colon \hat{\mcC^{op}}\to (Ab)$ is exact. 
We denote by $ {\rm Flat} (\mcC^{op},Ab)$ the full subcategory of flat functors. 
\end{dfn}

\begin{thm} (\cite{CB}, Thm p. 1646)\label{CBflat}\\
$\overrightarrow{\mcC}={\rm Flat} (\mcC^{op},Ab)$ and $F\in \overrightarrow{\mcC}$ is equivalent to 
\begin{itemize}
\item[(1)] $\mcC/F$ is filtered
\item[(2)] Every natural transformation $\coKer (-,f) \to F$ factors over a representable.  
\end{itemize}
\end{thm}

\begin{dfn}
Let $\mcA$ be an additive category. We say an object $X$ in $\mcA$ is \textbf{finitely presented} if $\Hom_{\mcA} (X,-)$ commutes with arbitrary filtered colimits. We denote by ${\rm fp} (\mcA)$ the full subcategory of finitely presented objects in $\mcA$. \\
The additive category $\mcA$ is called 
\textbf{locally finited presented} if ${\rm fp} (\mcA)$ is essentially small and $\mcA$ is equivalent to $\overrightarrow{{\rm fp}(\mcA)}$.  
\end{dfn}

\begin{rem}
If $\mcA$ is locally finitely presented then ${\rm fp} (\mcA)$ is essentially small, closed under direct sums and summands. In particular by Lemma \ref{ic}, it is idempotent complete.
\end{rem}

\begin{lem} (cf. \cite{CB}, part of Thm on p.1647) \\
For an essentially small category $\mcC$ we have ${\rm fp} (\overrightarrow{\mcC})\cong \mcC^{ic}$ is equivalent to the idempotent completion of $\mcC$.   
\end{lem}

\begin{thm} (\cite{CB}, Thm. in (1.2), p.1645)
If $\mcC$ is essentially small, then 
\[
{\rm fp} (\hat{\mcC})=\{F\in \hat{\mcC} \mid F\cong \coKer (-,f), \; f \in {\rm{Mor}} (\mcC)\}\quad  =:\modu_1\mcC
\]
Furthermore, $\hat{\mcC}$ is locally finitely presented. 
\end{thm}

 \begin{exa} Locally finitely presented abelian categories are \textbf{Grothendieck categories} (i.e. (1) abelian, (2) with arbitrary small coproducts, (3) directed colimits are exact, (4) has a generating object $G$). Here, the generator can be chosen as 
 \[
 G=\bigoplus_{C\in {\rm Ob}(\mcC)} (-,C)  \quad \in \overrightarrow{\mcC} 
 \]
 For the converse: If a Grothendieck category admits a set of finitely presented objects whose coproduct is a generator, then it is locally finitely presented. \\
 Grothendieck categories always have enough injectives (often they are hard to find), have arbitrary small limits and colimits. 
 \end{exa}

\begin{rem}
For a not necessarily small category $\mcC$ we can still define its ind-completion. If $\mcC$ is abelian, this is an abelian category - but it may not have enough injectives (cf. \cite{KS-Ind-Sheaves}, Prop. 15.1.2). 
\end{rem}

The following is a consequence of Lem. \ref{UP-add} together with \cite{CB}, Thm in (1.4), p. 1647. 
\begin{thm} (equivalence of (2-)categories) \\
The assignments $\mcC\mapsto \overrightarrow{\mcC}$ and $\mcA\mapsto {\rm fp} (\mcA)$ are $2-$functorial and give rise to an equivalence of  (strict) $2-$categories between 
\begin{itemize}
\item[(1)] essentially small, idempotent complete additive categories $\mcC$ with additive functors
\item[(2)] Locally finitely presented additive categories $\mcA$ with
additive functors that preserve arbitrary filtered colimits and restrict to the subcategories of finitely presented functors. 
\end{itemize}
\end{thm}

Let $\mcC$ be idempotent complete, essentially small additive category and $\mcA$ a locally finitely presented (additive) category. 
We assume $\mcC={\rm fp} (\mcA)$ and $\mcA= \overrightarrow{\mcC}$. 
Then the following holds (by restricting further and further): 
\begin{itemize}
\item[(i)] $\mcC$ left abelian $ \Leftrightarrow$ $\mcA$ abelian\\
the definition of \textbf{left abelian} (cf. \cite{CB}, (2.4)): Every morphisms has a cokernel, every epi is a cokernel and whenever $A\xrightarrow{f} B \xrightarrow{c} C$ with $c=\coKer (f)$ and $g\colon D\to B$, $cg=0$ then there exists an epi $d\colon E\to D$ such that $gd$ factors over $f$. 
\item[(ii)] $\mcC$ abelian $ \Leftrightarrow$ $\mcA$ locally coherent
\item[(iii)] $\mcC$ abelian and all objects noetherian $\Leftrightarrow$ $\mcA$ locally noetherian abelian
\item[(iv)] $\mcC$ is length abelian $\Leftrightarrow$ $\mcA$ is locally finite abelian
\end{itemize}

\begin{exa}
$R$ a ring, $R\Modu$ the category of all left $R$-modules and $R\modu_1$ the catego\\
(i) $R\Modu $ abelian and $R\modu _1$ is left abelian\\
(ii) for $R$ left coherent\\
(iii) for $R$ left noetherian\\
(iv) e.g. for $R$ left artinian (with Loewy length)
\end{exa}


\subsection{Gabriel Quillen embedding}

We review this well-known embedding of an essentially small exact category as a fully exact subcategory in an abelian category. \\

Let $\mcE=(\mcC, \mcS)$ be an essentially small exact category. 
\begin{dfn}
We define the category $\rm{Lex} (\mcE^{op}, Ab)$ to be the category of all additive functors $F\colon \mcC^{op} \to (Ab)$ which map short exact sequences $X \xrightarrow{i} Y \xrightarrow{d}Z$ in $\mcE$ to 
left exact sequences $0 \to F(Z) \xrightarrow{F(d)} F(Y) \xrightarrow{F(i)} F(X)$ in abelian groups. We will call this the category of \textbf{left exact functors on} $\mcE$.
We define the category of \textbf{locally effaceble functors} $\rm{Eff}_{\mcE}$ to be the full subcategory of $\hat{\mcC}$ of objects $F$ such that for every pair $(X, x)$ of an object $X$ in $\mcC$ and $x\in F(X)$ there exists an $\mcE$-deflation $d\colon Z \to X$ with $F(d) (x)=0$.  
\end{dfn}

\begin{lem} \cite{Kr-book}, Prop. 2.3.7 (1),(2) (with intermediate steps) 
\begin{itemize}
\item[(i)] (Prop 2.2.16) $\mcD=\rm{Eff}_{\mcE}$ is a Serre subcategory of $\hat{\mcC}$ closed under coproducts. Therefore, the 
Serre quotient functor $Q \colon \Modu \mcC \to \Modu \mcC/\mcD$ admits a right adjoint. 
\item[(ii)] (Lem. 2.2.10) Let $Q_{\rho}$ be the right adjoint. It factors as $\hat{\mcC} /\mcD \xrightarrow{\Phi} \mcD^{\perp} \xrightarrow{I} \hat{\mcC}$ with $\Phi $ an equivalence of categories and $I$ the inclusion functor. The quasi-inverse of $\Phi$ is given by $Q\circ I$. 
\item[(iii)] $\rm{Lex} (\mcE^{op}, Ab)=\mcD^{\perp}:=\{Y \in \hat{\mcC} \mid \Hom_{\Modu \mcC} (E,Y)=0= \Ext^1_{\Modu \mcC}(E,Y) \forall E\in \mcD\}$. 
\end{itemize}
\end{lem}
\begin{rem}
 $\rm{Lex} (\mcE^{op}, Ab)$ is a Grothendieck category (as it is the localization of a Grothendieck category by a Serre subcategory). In the abelian structure it has a generator 
 \[
 G=\bigoplus_{X\in {\rm Ob} (\mcC) } (-,X)
 \]
 As $(-,X)$ are (some) finitely presented objects in ${\rm Lex}(\mcE^{op},Ab)$, it follows that ${\rm Lex}(\mcE^{op},Ab)$ is locally finitely presented abelian. \\ 
 It also has an exact substructure as fully exact category in $\hat{\mcC}$ but these two exact structures usually do not coincide. 
\end{rem}


\begin{rem}
The inclusion $\rm{Lex} (\mcE^{op}, Ab)\to \hat{\mcC}$ is not an exact functor (if we consider $\rm{Lex} (\mcE^{op}, Ab)$ with its abelian structure). Yet it reflects exactness in the following sense: If $0 \to F \to G \to H \to 0$ in $\hat{\mcC}$ with $F,G,H$ in $\rm{Lex} (\mcE^{op}, Ab)$, then this is a short exact sequence in $\rm{Lex} (\mcE^{op}, Ab)$.    
\end{rem}

\begin{cor}
$\rm{Lex} (\mcE^{op},Ab) \subseteq \hat{\mcC}$ is a deflation-closed subcategory: 
Given a short exact sequence $0 \to F \to G \to H \to 0$ in $\Modu \mcA$. \\
If $G,H\in \rm{Lex} (\mcE^{op}, Ab)$, then also $F \in \rm{Lex} (\mcE^{op}, Ab)$. \\
If $F,G \in \rm{Lex} (\mcE^{op}, Ab)$ and $\Ext^2(E,F)=0$ for all $E$ effaceable, then $H \in \rm{Lex} (\mcE^{op}, Ab)$.  
\end{cor}

In general we can characterize short exact sequences in $\rm{Lex}(\mcE^{op},Ab) $ (in the abelian structure) as follows:  
\begin{lem}
Given two composable maps $0 \to F \xrightarrow{i} G \xrightarrow{p} H \to 0$ in $\rm{Lex}(\mcE^{op},Ab) $. TFAE
\begin{itemize}
\item[(1)] $(i,p)$ are an exact sequence in $\rm{Lex}(\mcE^{op},Ab) $  
\item[(2)] In $\hat{\mcC}$ we have a exact sequence $0 \to F \xrightarrow{i} G \xrightarrow{p} H $ and $\coKer (p) $ is locally effaceable.  
\end{itemize}
\end{lem}

\begin{thm} Let $\mcE=(\mcC, \mcS)$ be an exact category. 
The Yoneda functor gives a functor 
\[
i\colon \mcE \to \rm{Lex} (\mcE^{op}, Ab), \quad X\mapsto  (-,X)=\Hom_{\mcC} (-,X)
\]
with the following properties
\begin{itemize}
\item[(1)] $i$ is exact, reflects exactness and the essential image of $i$ is extension-closed, its idempotent completion is deflation-closed.  
\item[(2)] $i$ induces isomorphisms on all extension-groups. \\
\end{itemize}
\end{thm}

\begin{proof}
(1) \cite{Bue}, Appendix A and \cite{Kr-book}, Prop. 2.3.7.(3). The last statement follows from \cite{P-locCoh}, Prop. 6.1, (e).\\   
(2) \cite{Kr-book} in Lem 4.2.17 it is shown that $\mcE \to \rm{Lex} (\mcE^{op}, Ab)$ is right cofinal (Keller's definition) that implies the statement. 
\end{proof}

The following result is also relevant for us: 
\begin{pro} (\cite{P-locCoh}, Prop. 6.2) \label{PosExt1}
For every $X$ in the category $\mcE$, the functor \\$\Ext^n_{{\rm Lex} (\mcE^{op}, Ab) }((-,X), - )$ preserves filtered colimits. 
\end{pro}

\subsection{Locally coherent exact}

The ind-completion of a small exact category $\mcE$ has a natural exact structure (namely as fully exact in the Gabriel-Quillen embedding), this exact structure can also be described as directed colimits of short exact sequences in $\mcE$ and is called \textbf{locally coherent exact structure}.\\ 
The CB-correspondence from the previous subsection is 'extended' to i.c. small exact categories.

We begin with the following observation. 
\begin{lem}
Let $\mcC$ be a small additive category, $F$ a flat functor and  $X\xrightarrow{f} Y \xrightarrow{g} Z$ with $g=\coKer (f)$ in $\mcC$ then $0 \to F(Z) \to F(Y) \to F(X)$ is exact in abelian groups.     
\end{lem}

\begin{proof}
Define $(X, -):=\Hom_{\mcC} (X,-)\colon \mcC \to (Ab)$ and the contravariant Yoneda embedding $\mcC^{op}\to \hat{\mcC^{op}}, X\mapsto (X,-)$. By assumption we have an exact sequence $0 \to (Z,-) \to (Y,-) \to (X, -)$. Since $F$ is flat, the functor $F\otimes_{\mcA}$ is exact and we have $F\otimes_{\mcC} (E,-)\cong F(E)$. Therefore, we obtain the exact sequence $0 \to F(Z) \to F(Y) \to F(X)$.     
\end{proof}

Now for an exact category $\mcE=(\mcC, \mcS)$, by the previous lemma $\rm{Flat} (\mcC^{op},Ab)\subseteq {\rm Lex}(\mcE^{op},Ab)$.  
\begin{lem} (and definition.)\label{IndE}
   $\rm{Flat} (\mcC^{op},Ab)$ is closed under extensions in ${\rm Lex}(\mcE^{op},Ab)$. \\
   We define $\overrightarrow{\mcE}$ to be the fully exact structure on $\overrightarrow{\mcC}$ and call this the \textbf{ind-completion of the exact category} $\mcE$.   
\end{lem}

\begin{proof} (sketch) (of Lemma \ref{IndE})
Given a short exact sequence $\sigma \colon F \to G \to H$ in (the abelian structure on) ${\rm Lex}(\mcE^{op},Ab)$ with $F,H$ flat.  \\
First assume $H=(-,X)$, write $F$ as filtered colimit and use Prop \ref{PosExt1} and the fact that the essential image of $\mcE\to {\rm Lex}(\mcE^{op},Ab)$ is extension-closed to conclude the claim. \\
In general, we use Thm \ref{CBflat}, (2). Given a morphism $\theta\colon \coKer (-,f) \to G$. Postcompose to $\coKer (-,f)\to H$. As $H$ is flat, this factors over a morphism $g\colon (-,X)\to H$ for some $X$ in ${\rm Ob}(\mcE)$. Now form the pull-back of $\sigma$ along $g$ in the abelian category ${\rm Lex}(\mcE^{op},Ab)$, say this is a short exact sequence $F\to E \to (-,X)$. 
The universal property of the pull-back gives a morphism $\theta'\colon \coKer (-,f) \to E$ and a morphism $u \colon E \to G$ with $\theta = u\theta'$. By the first case $\theta'$ factors over a representable, therefore $\theta $ does so too.  
\end{proof}

\begin{rem} (and definition)
Let $\mcE=(\mcC, \mcS)$ be an essentially small exact category, then $\mcE$ is fully exact in ${\rm Lex}(\mcE^{op},Ab)$. As the latter is abelian, it is idempotent complete, therefore $\mcE^{ic}$ is also a fully exact subcategory in ${\rm Lex}(\mcE^{op},Ab)$.\\ 
But this means $\rm{fp}(\overrightarrow{\mcC})$ is an extension-closed subcategory in $\overrightarrow{\mcE}$. We define $\rm{fp}(\overrightarrow{\mcE})$ to be the fully exact structure on $\rm{fp}(\overrightarrow{\mcC})$. 
\end{rem}

\begin{dfn}
Let $\mcF=(\mcA, \mcT)$ be an exact category. We say it is \textbf{locally coherent exact} if $\rm{fp}(\mcA)$ is essentially small and extension-closed in $\mcF$ -in this case we denote by ${\rm fp} (\mcF)$ the fully exact subcategory- and $\mcF=\overrightarrow{{\rm fp} (\mcF)}$.
\end{dfn}

In an ess. small exact category $\mcE$, the category of short exact sequences ${\rm Ses}(\mcE)$ is again an essentially small exact category (with degree-wise short exact sequences). 
\begin{thm} \label{PosThm} (\cite{P-locCoh}, proof of Lemma 1.2) 
A filtered colimit of short exact sequences in $\overrightarrow{\mcE}$ is a short exact sequence in $\overrightarrow{\mcE}$.\\
The universal property of the ind-completion yields an equivalence of categories
\[
\overrightarrow{{\rm Ses} (\mcE)} \to {\rm Ses} (\overrightarrow{\mcE})
\]    
\end{thm}

We also observe the following: 
\begin{lem}
For an essentially small exact category $\mcE$, we have that $\overrightarrow{\mcE}$ is a resolving subcategory in ${\rm Lex}(\mcE^{op},Ab)$. In particular it is homologically exact. 
\end{lem}

\begin{proof}
$\overrightarrow{\mcE}$ is extension-closed, idempotent complete and contains a generator of ${\rm Lex}(\mcE^{op},Ab)$ , so it is enough to see that it is also deflation-closed. Now given a short exact sequence $F \to G \to H$ in ${\rm Lex}(\mcE^{op},Ab)$, it gives a $4$-term exact sequence in $\hat{\mcC}$: $0\to F\to G \to H \to E\to 0$ withe $E$ locally effaceable. Assume that $G,H$ are flat, we want to see that $F$ is so too. As every finitely presented effaceable functor has projective dimension $2$, every locally effaceable has flat dimension $2$ and this implies $F$ is flat.
\end{proof}

Furthermore, since $(-,X)\in \rm{Flat} (\mcE^{op},Ab)$ for all $X$ in $\mcC$, we get a fully faithful exact functor with extension-closed essential image 
\[
i\colon \mcE \to \overrightarrow{\mcE}, \quad X \mapsto (-,X)
\]

\begin{lem}
The functor $i$ is homologically exact. 
\end{lem}
\begin{proof}
As $\mcE \to {\rm Lex}(\mcE^{op},Ab)$ is homologically exact and also $\overrightarrow{\mcE}\to {\rm Lex}(\mcE^{op},Ab)$, this is immediate. 
\end{proof}

We directly get the following from the previous corollary and Prop. \ref{PosExt1}. 
\begin{cor} 
For $X$ in ${\rm Ob} (\mcE)$, the functor $\Ext^n_{\overrightarrow{\mcE}}((-,X), -)$ commutes with filtered colimits.     
\end{cor}

Let ${\rm Ex}(\mcE, \mcF)$ be the category of exact functors between two exact categories $\mcE$ and $\mcF$ and ${\rm Ex}_{fc} (\mcE, \mcF)$ be the full subcategory of exact functors which preserve filtered colimits. 
\begin{lem} (\textbf{Universal property of ind-completion for exact categories}) \label{UElocCoh}
Let $\mcE$ be an essentially small exact category. 
Then $\overrightarrow{\mcE}$ is closed under directed colimits and directed colimits are exact functors. \\
Let $\mcF$ be an exact category closed under all directed colimits and they are exact functors. Then precomposition $\mcE\to \overrightarrow{\mcE}$ gives an equivalence
\[{\rm Ex}_{fc} (\overrightarrow{\mcE}, \mcF)\to {\rm Ex} (\mcE, \mcF)\] 
\end{lem}

\begin{thm} (equivalence of 2-categories) \\
The assignments $\mcE\mapsto \overrightarrow{\mcE}$ and $\mcF\mapsto {\rm fp} (\mcF)$ are functorial and give rise to an equivalence of $(2-)$categories between 
\begin{itemize}
\item[(1)] essentially small, idempotent complete exact categories $\mcE$ with exact functors 
\item[(2)] Locally coherent exact categories $\mcF$ with
exact functors that preserve arbitrary filtered colimits and restrict to the subcategories of finitely presented functors. 
\end{itemize}
\end{thm}
We remark that a functor that preserves filtered colimits on objects, also preserves filtered colimits on morphism categories and this implies it also preserves filtered colimits of short exact sequences. This implies that as exact functor between essentially small exact categories $F\colon \mcE_1\to \mcE_2$ extends with the universal property of the ind-completion uniquely to an exact functor $\overrightarrow{F}\colon \overrightarrow{\mcE_1}\to \overrightarrow{\mcE_2}$, this is a consequence of Thm \ref{PosThm}. 
Recall: $F$ fully faithful if and only $\overrightarrow{F}$ is fully faithful by Prop. \ref{SGA4Ind}. 

Ignoring set-theory for a moment: Let $\mcF$ be an exact category, 
we consider ${\rm EX} (\mcF)$ the lattice of all exact subcategories. 
For $\mcF$ locally coherent exact, 
we define ${\rm EX}_{fc} (\mcF)$ to be the exact subcategories in the category (2) above, i.e. exact functors $i \colon \mcF'\to \mcF $ such that $i$ is fully faithful, $\mcF'$ is locally coherent exact and $i$ preserves filtered colimits.   
\begin{cor} Let $\mcE$ be an essentially small category. We have mutually inverse, isomorphisms of posets
\[
\overrightarrow{(-)} \quad \colon {\rm EX} (\mcE) \leftrightarrow {\rm EX}_{fc} (\overrightarrow{\mcE})\colon\quad  {\rm fp}(-)
\]
It restricts to all the usual subposets such as extension-closed, exact substructures etc. 
\end{cor}

Positselski found the maximal and minimal locally coherent exact structure on a locally finitely presented category particular interesting. The minimal exact structure is the ind-completion of the split exact structure (on an essentially small additive category) and is called \textbf{pure exact structure} on a locally finitely presented category. 

\begin{exa}
Let $R$ be a ring. 
The abelian exact structure on $R\Modu $ is the maximal locally coherent exact structure, it corresponds to the left abelian structure on $R\modu _1$ (fp $R$-modules). 
The category of flat $R$-modules $R\Modu _{fl}$ is extension-closed in $R\Modu $. Its subcategory of finitely presented objects is $(\add(R))^{ic}$ with the split exact structure is the fully exact substructure. By \cite[Thm 3]{CB} $\overrightarrow{\add(R)}=R\Modu _{fl}$. In this case: The fully exact structure is the pure exact structure. 
\end{exa}

Stovicek generalized the notion of a Grothendieck category to an exact category of Grothendieck type. 
\begin{thm} (\cite[Cor. 5.4]{P-locCoh}) Locally coherent exact categories are exact categories of Grothendieck type (in the sense of Stovicek).  
\end{thm}
In particular, all established properties of exact categories of Grothendieck type hold true. 
\begin{cor} (also \cite[Cor. 5.4]{P-locCoh})$\overrightarrow{\mcE}$ has enough injectives. 
\end{cor}

\bibliographystyle{alpha}
\bibliography{Sauter-ZieglerExCats}
\end{document}